\documentclass[ece]{puthesis}

\usepackage{amsmath, amssymb, latexsym}
\usepackage{psfrag}

\usepackage{graphicx}
\usepackage{epsfig}
  \ifx\LabelFigloaded\MYundefined\relax
  \else
   \fi

  \def\LabelFigloaded{\relax}

  \chardef\LabelFigCatAt\the\catcode`\@
  \catcode`\@=11

 \let\LabelFigwlog@ld\wlog
 \def\wlog#1{\relax}

 \ifx\\\MYundefined@
    \let\\\relax
 \fi

  \def\ms@g{\immediate\write16}

 \def\N@wif{\csname newif\endcsname }
 \def\Temp@ {\N@wif\ifIN@}
 \ifx\INN@\MYundefined@
    \else \let\Temp@\relax
 \fi
 \Temp@

  \def\IN@{\expandafter\INN@\expandafter}
  \long\def\INN@0#1@#2@{\long\def\NI@##1#1##2##3\ENDNI@
    {\ifx\m@rker##2\IN@false\else\IN@true\fi}%
     \expandafter\NI@#2@@#1\m@rker\ENDNI@}
  \def\m@rker{\m@@rker}
 
  \newtoks\Initialtoks@  \newtoks\Terminaltoks@
  \def\SPLIT@{\expandafter\SPLITT@\expandafter}
  \def\SPLITT@0#1@#2@{\def\TTILPS@##1#1##2@{%
     \Initialtoks@{##1}\Terminaltoks@{##2}}\expandafter\TTILPS@#2@}

 \def\Shifted@@#1#2#3{\setbox0=\hbox{#3}%
   \raise -\dp0\vbox {\kern-#2%
       \hbox {\kern#1\unhbox0\kern-#1}%
           \kern#2}}

 \newcount\gridcount
 \newbox\auxGridbox@ \newbox\hGridbox@ \newbox\vGridbox@
 \newbox\Labelbox@ \newbox\auxLabelbox@
 \newbox\Coordinatebox@
 \newtoks\Labeltoks@
 \newdimen\Wdd@ \newdimen\Htt@
 \newdimen\Wddd@ \newdimen\Httt@
 
 \def\Wr@{\immediate\write16}

 \newdimen\GL@wd
 \def\GridLineWidth#1{\GL@wd=#1}

 \def\gobble#1{}
 \def\EdgeErr@{\Wr@{}%
      \Wr@{\string\Edges\space argument
      1, 10, 100 or 1000 please\string!}%
      }

 \newcount\Edgect@

 \def\Sweepup#1\endSweepup{}

 \def\SetEdges@{%
    \edef\Zr@@s{\expandafter\gobble\number\Edgect@\empty}%
        \count255=0\Zr@@s\relax
        \ifnum\count255=\z@\else\EdgeErr@\show\tailtest\fi
        \count255=1\Zr@@s\relax
        \ifnum\count255=\Edgect@\relax\else\EdgeErr@\show\leadtest\fi
    \EdgGl@b\edef\Zr@s{\expandafter\gobble\Zr@@s\empty}
    \ifnum\Edgect@>\@ne\relax\EdgGl@b\let\L@Dc\empty
        \else\EdgGl@b\edef\L@Dc{\string.}\fi
    \ifnum\Edgect@>\@ne\relax
        \EdgGl@b\edef\Edgescale@##1{\divide##1 by \Edgect@}%
        \else\EdgGl@b\edef\Edgescale@##1{}\fi
    }

 \def\Edges#1{\Edgect@=#1\relax
     \let\EdgGl@b\global \SetEdges@}

 \Edges{1}

 \def\hhrule{\hrule height \GL@wd\vskip-.\GL@wd}

 \def\hRule@{%
   \advance\gridcount -2%
   \vfil\hhrule\vfil
   \llap{\smash{\raise -2.5pt
     \hbox{\L@Dc\number\gridcount\Zr@s\kern2pt}}}%
   \hhrule
   }

\def\vvrule{\vrule width \GL@wd \kern-\GL@wd}

 \def\vRule@{\advance\gridcount 2%
   \hfil\vvrule\hfil
   \setbox\auxGridbox@=\vbox to 0pt
      {\vskip \Htt@\vskip 2pt
        \hbox to 0pt{\hss\L@Dc\number\gridcount\Zr@s\hss}\vss}%
      \wd\auxGridbox@=0pt \box\auxGridbox@
   \vvrule
   }

 \def\PlaceGrid@@{\gridcount=10 
  \setbox\hGridbox@=\hbox{%
        \hbox{%
             \hskip-.4pt\vrule
             \vbox to \Htt@{%
               \offinterlineskip\parindent=\z@\relax
               \hbox to \Wdd@{\hfil}
               \hRule@\hRule@\hRule@\hRule@
               \vfil\hhrule\vfil}%
             \vrule\hskip-.4pt}
    }%
  \gridcount=0%
  \setbox\vGridbox@=\hbox{%
      \vbox{\offinterlineskip\parindent=0pt\hsize=0pt
         \vskip-.4pt\hrule%
         \hbox to \Wdd@{%
                 \vtop to \Htt@{\vfil}%
                 \vRule@\vRule@\vRule@\vRule@
                 \hfil\vvrule\hfil}%
         \hrule\vskip-.4pt}}%
  \wd\hGridbox@=0pt\ht\hGridbox@=0pt
  \wd\vGridbox@=0pt\ht\vGridbox@=0pt
  \hbox{\box\hGridbox@\box\vGridbox@}%
  }

 \def\LabelsGlobal{\def\LabGl@b{\global}}
 \def\LabelsLocal{\def\LabGl@b{}}
 \LabelsGlobal 

 \def\SetLabels#1\endSetLabels{%
   \LabGl@b\Labeltoks@={#1()\\}%
   }

 \LabGl@b\Labeltoks@={()\\}

 \def\ShowGrid{\LabGl@b\let\PlaceGrid@\PlaceGrid@@}
 \def\HideGrid{\LabGl@b\let\PlaceGrid@\relax}
 \def\Grids{\ShowGrid\LabGl@b\let\GridSwitch@\ShowGrid}
 \def\noGrids{\HideGrid\LabGl@b\let\GridSwitch@\HideGrid}

 \noGrids

 \def\bAdjust@@{%
     \setbox\auxLabelbox@=\hbox{\raise \dp\auxLabelbox@
            \box\auxLabelbox@}}
 \def\bAdjust@{\let\vAdjust@\bAdjust@@}

 \def\eAdjust@@{\dimen0=-.5\ht\auxLabelbox@
     \advance\dimen0 by .5\dp\auxLabelbox@
     \setbox\auxLabelbox@=
            \hbox{\raise\dimen0\box\auxLabelbox@}}
 \def\eAdjust@{\let\vAdjust@\eAdjust@@}

 \def\tAdjust@@{%
     \setbox\auxLabelbox@=\hbox{\raise-\ht\auxLabelbox@
            \box\auxLabelbox@}}
 \def\tAdjust@{\let\vAdjust@\tAdjust@@}

 \let\vAdjust@\relax

 \def\lAdjust@{\let\hAdjust@\rlap}
 \def\rAdjust@{\let\hAdjust@\llap}

 \let\hAdjust@\relax\let\vAdjust@\relax

 \def\FetchLabel@#1(#2)#3\\{%
     \IN@0#2@@\ifIN@
        \setbox0=\hbox{\ignorespaces#1#3\unskip}%
        \ifdim\wd0>0pt
           \ms@g{}%
           \ms@g{ !!! Bad label(s)? !!!}%
           \message{ #1(#2)#3}%
        \fi
        \def\LabelMole@##1\endFetchLabel@{%
            \IN@0()\\@##1@%
            \ifIN@\def\Temp@{\FetchLabel@##1\endFetchLabel@}%
            \else\def\Temp@{}%
            \fi
            \Temp@
           }%
     \else
       \ignorespaces#1\unskip
       \setbox\auxLabelbox@=%
         \hbox to 0pt{\hss\ignorespaces\hAdjust@
          {\ignorespaces#3\unskip}\hss}%
       \vAdjust@
       \let\hAdjust@\relax\let\vAdjust@\relax
       \AugmentLabelBox@@{#2}%
       \ht\Labelbox@=0pt\dp\Labelbox@=0pt
       \let\LabelMole@\FetchLabel@%
     \fi\LabelMole@}

 \newtoks\XYSep@ 
 \def\SetXYSeparator#1{%
     \IN@0#1@@\ifIN@\XYSep@{*}%
     \else
     \XYSep@{#1}%
     \fi
     }

 \SetXYSeparator*

 \def\AugmentLabelBox@@#1{%
     \IN@0\the\XYSep@ @#1@\ifIN@
       \SPLIT@0\the\XYSep@ @#1@%
       \setbox\Labelbox@=\hbox to 0pt{%
         \unhbox\Labelbox@
         \Shifted@@{\the\Initialtoks@\Wddd@}%
         {\the\Terminaltoks@\Httt@}%
         {\box\auxLabelbox@}}%
     \else
         \ms@g{}%
         \ms@g{ !!! Bad insertion point. !!!}%
         \message{ (#1\ this point was rejected.)}%
     \fi
    }

 \def\FetchOption@#1[#2]#3\endFetchOption@{%
    \def\temp{#1}
    \ifx\temp\empty
       \Edgect@=#2\relax
       \let\EdgGl@b\relax
       \SetEdges@
       \Cleaner@#3%
    \fi}

 \def\Cleaner@#1[@]{\Labeltoks@{#1}}
     
 \def\PlaceLabels@@{\mathsurround=0pt
     \def\Cr@{\\}%
     \let\L\lAdjust@\let\R\rAdjust@
     \let\B\bAdjust@\let\E\eAdjust@\let\T\tAdjust@
     \expandafter\FetchOption@\the\Labeltoks@[@]\endFetchOption@
     \Wddd@=\Wdd@ \Edgescale@\Wddd@ 
     \Httt@=\Htt@ \Edgescale@\Httt@
     \expandafter\FetchLabel@\the\Labeltoks@\endFetchLabel@
     \box\Labelbox@
     }%

 \let \PlaceLabels@\PlaceLabels@@

 \def\AffixLabels#1{\setbox\Coordinatebox@=\hbox{#1}%
      \Wdd@=\wd\Coordinatebox@ \Htt@=\ht\Coordinatebox@
      \advance\Htt@ \dp\Coordinatebox@
      \hbox{\copy\Coordinatebox@\kern-\Wdd@ 
           \Shifted@@{0pt}{-\dp\Coordinatebox@}%
           {\PlaceLabels@\PlaceGrid@}%
           \kern\Wdd@}%
      \GridSwitch@ 
      \LabGl@b\Labeltoks@{()\\}%
      }
 
   \let\wlog\LabelFigwlog@ld   
   \catcode`\@=\LabelFigCatAt  

                                By

              Raymond S\'eroul <A18645@FRCCSC21.BITNET>
                                and 
              Laurent Siebenmann <lcs@topo.math.u-psud.fr>
    
              VERSIONS: July 1991, Oct 1991, Jan 1992, July 1992

INTRODUCTION

      This labelling package is intended for TeX users who
rely on non-TeX sources for for their graphics inserts.  It
provides means for adding TeX labels to such inserts with a
minimum of fuss. 

       For most labels, TeX users have in the past found it
reasonably convenient to rely on non-TeX sources. Typical
occasions when an inescapable need for TeX labels seemed to
arise are

 (a) when the graphics program lacks certain exotic or complex
mathematical symbols

 (b) when the very highest typographical quality is wanted for the
labels

 (c) when labels included with the graphics fail to print, 
 and you cannot figure out why (cf. boxedeps.doc).  The labels
 provided by labelfig.tex are 100

       Since this package first appeared, many users, who in the
past scarcely dreamed of using TeX labels, have come to use
nothing but.  So it is now appropriate to add

Intoxication Warning:  TeX labels may be addictive and expensive. 

     If you have a fast preview you may disagree, and even find
that this package provides an agreeable paste-up environment; see
extra applications at end.

     Note to publishers: It is possible and convenient to ultimately
export the TeX labels produced by labelfig.tex to become an integral
part of the EPS file. This is often desired by a publisher who typically
uses an "upmarket" graphics or page layout program, with which the
staff is skilled in perfecting figures.  See Appendix I for
a recipe.

     The authors are grateful to Patrick Ion of Math Reviews for
helpful comments and encouragement.

BASIC INSTRUCTIONS

    After reading in the macro file using

preview or proof your figure with a coordinate grid printed on
top, by typing the following:

    \ShowGrid  
    \AffixLabels{<the graphics insertion>}

Here <the graphics insertion> is what you would type to insert
the graphics object alone without the grid.  This must provide
for the space around it. For example <the graphics insertion>
might well be \BoxedEPSF{MyFigure scaled 700} using the
boxedeps.tex macro package (from same source); this provides a
TeX box containing the encapsulated PostScript insert specified by
the file MyFigure. \AffixLabels{...} provides the grid (supposing
\ShowGrid is present) and later, once you have specified labels
using the grid, it will "tack on" the labels.

     The grid is a sort of (usually elongated) checkerboard of
ten rows and ten columns and its (internal) partitions are by
default numbered  .1, ... ,.9  both horizontally (X-coordinate
running left to right) and vertically (Y-coordinate running bottom
to top).  Thus the points enclosed by the grid correspond to the
points of the unit square in the cartesian "X-Y" plane, the lower
left corner corresponding to the origin (0,0).  By extrapolation,
the full page corresponds to a larger rectangle in the plane.

     These coordinates serve to position labels as follows.
Before the \AffixLabels{...} command type label specifications:

  \SetLabels
   (<X-coordinate>*<Y-coordinate>) <first label> \\
   .
   .
   .
   (<X-coordinate>*<Y-coordinate>)  <last label> \\
  \endSetLabels

Each row specifies one label and is terminated by \\.  In each
row, the position indicator comes first; it is written as a
standard cartesian point except that the X- and Y- coordinates
are separated by * rather than a comma because TeX allows a
comma as decimal point. There are no dimension units to specify
as the unit is the grid itself.

     By default, this cartesian point specifies where the middle
of the baseline of the label will be located.  However if you precede
the point by \L [or \R] the left [or right] edge of the baseline will
be located there. Similarly you may also precede the point by \T, \E,
or \B to vertically align the top equator or bottom of the label box
at the specified point.  This gives nine standard positions of
the label with respect to the insertion point --- corresponding to
the eight principle points of the compas and the center

                     \L\T     \T      \R\T

                     \L\E     \E      \R\E

                     \L\B     \B      \R\B

But this neglects the default "baseline" level of TeX,
giving potentially three more positions

                     \L    <no tag>   \R

For text, the baseline level is often the preferred. Its relation to
the others is variable. It will often coincide with the bottom level,
as happens for "X".  But it is often distinct, as for "g", in which
case you have in all 12 distinct positions rather than 9.

     It is convenient to think of this specification of label
position as attaching the label by a thumb-tack to the coordinate
grid. There are up to twelve positions of the thumb-tack on the
label, while the position of the thumb-tack on the coordinate grid is
arbitrary.  Normally, one choses the position of the thumb-tack on
the label to be the one that is the closest to the item being
labeled.  There are good reasons for this "rule of thumb":

   (a)  It facilitates correct positioning at first try.

   (b)  If the scale of the figure must be altered after labels
have been affixed, the labels have a good chance of remaining well
positioned.

   (c)  The visible grid need not extend beyond the "bounding box"
for the figure, because the best preferred position is always
(at least almost) within the bounding box .

The second reason is particularly important. Indeed it often
happens that scale has to be altered after labelling begins, in
order to either provide space for the labels, or to adjust
proportions between the labels and the figure.  (The size of labels
is unaffected by scaling.)

     Here is an artificial but self-contained test which uses
TeX rules to make a graphics object.

TEST

    Do not skip this!

 \def\FrameIt#1{\hbox{\vrule$\vcenter {\hrule\kern3pt%
             \hbox {\kern3pt #1\kern3pt}%
               \kern3pt\hrule}$\relax\vrule}}

 \def\Caption#1#2{\FrameIt{%
       \vtop {\hsize=#1\relax \parindent=0pt
         \leftskip=0pt \rightskip=0pt plus15pt
         \parfillskip=0pt
         \lineskip=1pt\baselineskip=0pt
         #2}}}

 \def\FirstQuadrant{\hbox to 100pt{\vrule\vbox to 100pt{%
        \hbox to 100pt{\hfil}\vfil\hrule}\hss}}

  \SetLabels
    \R(.5*.2) $\zeta\,\cdot$\\
    (.9*-.10) $\xi$\\
    \R(-.03*.9) $\eta$\\
    \T(.5*.9) \Caption{70pt}{%
          \it The norm of
          $g(\xi+i\eta)$ is indicated on
          contours of this invisible surface.}\\
  \endSetLabels

  \AffixLabels{\FirstQuadrant} \end

  Note that the coordinates to use for labels are indicated on the
edges of the grid (when visible) corresponding to the conventional
x- and y- axes of the Cartesian plane. By default the grid is
1-by-1. However, by the command \Edges{100}, you can change this
to 100-by-100 and many users find this alternative most
convenient. Place the command \Edges{...} in your style file (or
header) since its effect is is global. Other possible edge values
are 10 and 1000.

  If you use the command \Edges{...} at all, do so with care.  For
if you accidentally delete an \Edges{...} command your labels will
abruptly be badly misplaced and may logically but mysteriously
generate "dimension too big" errors under TeX and "off page" errors
under your driver.  

  You can dictate the edgescale for an individual figure by giving
the scale in brackets immediately after \SetLabels.  Thus, to
import into an article using say \Edge{100} a figure labelled using
another edgescale, say the original 1-by-1 default, you can use
\SetLabels[1]...\endSetLabels.

GETTING IT DOWN PAT

     Complicated labeling deserves the same respect as
complicated mathematics.  Do not expect it to come out perfect the
first time!  What is needed in either case is a mechanism to
repeatedly typeset troublesome pieces.

     One mechanism is always available.  One does complicated
labelling in a separate "test" file involving just the figure being
labelled;  a texpert will know how to \dump TeX's current state as
a temporary format that restarts rapidly at each retry.  Usually,
one then pastes the completed labelled figure back into the main
TeX file, but, of course, one can also \input it as an auxiliary
file.

     If you do not have a TeXpert at handy, here is a first
approximation to an efficient setup. By deletions reduce a copy
of your article to just a few lines before and after the figure.
Now label the figure, and finally, copy and paste the labelled
figure to the original article. Then copy the next figure to label
into this testbed and repeat. The TeXpert can improve the  speed
at which TeX starts up, by compiling a format specifically for
your article; just one caution: best NOT include in the format
ephemeral details of setup like \Set<mydriver>ArtSpecials (from
boxedeps.tex because this reads  figure dimensions which you may
change during your work session.

     An improved mechanism to repeatedly typeset troublesome
pieces is now available on the Macintosh; it is called LinoTeX;
see the same ftp sources.  It could be set up on many types
of computer.

     Before using labelfig.tex to attach labels to a graphics
object inserted using boxedeps.tex or BoxedArt.tex, make it a
firm rule to carefully adjust the bounding box using the trimming
commands of these packages, and also at least tentatively scale
and position the object. Beware of changing the grid inadvertently
after the labels have been positioned.  For example, correcting
the bounding box of a PostScript graphics object can foul up the
labels by changing the coordinate grid to which the labels are
attached. This is particularly true for the trimming  commands of
boxedeps.tex and BoxedArt.tex. However, as noted already, change
of scale is much less disruptive, and modest adjustments should be
well tolerated.

     Sometimes the labels protrude so far from the bounding box
of a figure that the figure has to be repositioned.  Best do this
by ad hoc spacing, say using \hglue and \vglue; altering the
bounding box would create a vicious circle.

     Remember that you are responsible for preventing labels
from overlapping. You are responsible for all label typography
including size and style. A label is really just about anything
that can be put in a TeX box. Note that spaces at the beginning
and end of labels will normally be suppressed; if you really want
them you must protect them with TeX braces.

     This package temporarily sets the \mathsurround parameter
of TeX to zero  while the labels are being affixed. This is done
because nonzero \mathsurround space would influence the position
of left and right aligned labels; then, when a texpert or printer
modifies mathsurround, diagram labeling might be disastrously
altered. There is a small price to pay involving labels that are
formatted as caption boxes including mathematics: you  may want or
need to specify an explicit mathsurround space within the caption
box; it will not influence anything outside.

     Those hostile to the use of * as separator between
the X and Y coordinates of label insertion points, are free to
impose another using \SetXYSeparator{<the new separator>}.  
Americans may prefer "," to "*" since they never use a 
comma as a decimal point; on the other hand, * may be more visible.

APPENDIX (I)  MERGING labelfig.tex LABELS INTO AN EPSF GRAPHICS OBJECT.

     As promised in the introduction, here is a recipe useful for
publishers. It works at least on Macintosh and at least for vectorized
graphics and Adobe type1 fonts.  (There is surely a similar recipe for
PCs under MSWindows.)

 (a)  Use boxedeps.tex utility to integrate the figure given by the eps
file, "x.eps" say, with a visible frame around it.  See
\ShowDisplacementBoxes command in boxedeps.tex.  To get precise results
automatically it is important to use the \Trim... commands of
boxedeps.tex making the "DisplacementBox" neatly fit the figure.

 (b)  Use the TeX printer driver and LaserWriter (versions >= 8.1.1) to
export to an EPSF the DVI page containing the integrated, labelled
figure. You now have an EPS file  "xx.eps"  that contains too much, and at
the wrong scale, and at wrong position.

 (c)  Convert the EPSF to an Adode Illustrator format EPSF using
the shareware utility called epsConvert by Sam Weiss
1993-- (currently $25).

 (d)  In Illustrator (or a compatible program), group the labels and the
"DisplacementBox"; copy them to the clipboard and paste them into "x.ps".
This step requires that all the label fonts be "visible to the Macintosh.

 (e)  Translate and scale the pasted group consisting of the labels plus
the "DisplacementBox" so as to make the "DisplacementBox" the bounding
box of (labelless) figure represented by "x.eps".  At this point the
labels will be correctly placed on the figure "x.eps".

 (f)  Ungroup and delete the "DisplacementBox".  The result is the
desired single EPS file, "x+.eps" say, It contains the original figure
plus its labels.  

     Using grouping and ungrouping appropriately in "x+.eps", a
publisher's staff can very efficiently improve label positions etc.

APPENDIX II)  SOME EXOTIC APPLICATIONS

     The grid of labelfig.tex is analogous to a light-table in
classical page makeup with wax or latex glue.  In principle, you
can use it to compose any page from its indivisible parts.  This
even has some of the artisanal charm of classical paste-up
provided you have a fast screen preview to make the process
"interactive".

     In practice labelfig.tex is a tool for nonstandard jobs.
Here are a few going beyond the labelling already discussed.

(I)  GRAPHICS INTEGRATION.

     This is accomplished by treating the imported graphics
objects as labels.  The underlying graphics object is then
typically an empty  \vbox to <dimension>{\vfill} in a TeX
\midinsert...\endinsert construction.  A label line
might be of the form

   (.1*.1) \\

The exact form of the special command varies from driver to
driver.  However, in the case of encapsulated PostScript graphics
(EPSF norm), by relying on boxedeps.tex, one can have the
following standard syntax (independant of driver  (see
boxedeps.doc for details.
  
  (.1*.1) \BoxedEPSF{MyFigure scaled <scale in mils>}\\

This may be slow since it requires TeX to read the PostScript
file to read bounding box using many complex macros.  So you
may want to try

  (.1*.1) \EPSFSpecial{MyFigure}{<scale in mils>}\\

which is fast and driver independant, but it squashes the
bounding box, normally to its lower left corner.

     Similarly for graphics of the Macintosh PICT norm ---
using BoxedArt.tex (same sources) in place of boxedeps.tex.

     This approach to integration is to be recommended when
one is assembling a composite graphics object.

 (II)  COMMUTATIVE DIAGRAM ENHANCEMENT

     Commutative diagrams or arrays of mathematical objects
connected by arrows of various sorts are common in mathematics.
The mathematical objects require the use of TeX.  Recently TeX
acquired a good collection of arrows of all slopes --- that of
LamSTeX --- plus pwerful macros to build the diagrams.

     However, even the LamSTeX collection is often
inadequate; it lacks for example double shafted arrows, dotted
arrows and curved arrows. Fortunately it is possible to produce
such arrows on an individual basis using sophisticated graphics
programs such as Illustrator and AldusFreehand (both serving
the EPSF norm) or using Metafont (with its public domain norm).
Since the creation of each new arrow is a work of love, you
probably want to limit the number of arrows by using LamSTeX
for most arrows. The 40K commutative diagram module of LamSTeX
has been adapted to work with AmSTeX and a copy may be posted
with LabelFig and related files. Unfortunately no one has yet
offered a version that works with Plain TeX or LaTeX.

       Suffice it here to say that when the exotic arrow has
been somehow imported into TeX, labelfig.tex treats it as a
label that one affixes to the commutative diagram.  Two other
steps will be treated in separate notes, namely the matter of
extracting the dimension specifications for the arrow and the
construction of the arrow --- for these steps are far from
unique and often depend intimately on your computer environment. 
Notes for the Macintosh-Textures-Illustrator combination are
found in the file ExoticArrows.doc.

 (III) NESTING 

Ingenuity pays off in exploiting labelfig.tex. One can
mix graphics and typography quite freely.  labelfig.tex is good
for freeform or overlapping arrangements, while boxedeps.tex (or
BoxedArt.tex) is best for regimented non-overlapping
arrangements --- and the two can be combined.

     The default behavior of labelfig.tex is not ideal 
for nesting objects, because to prevent trouble for beginners
the register for labels is globally cleared when \AffixLabels
concludes.  But there are switches available

      \LabelsGlobal      \LabelsLocal

which change this.  To understand this, extend the above test 
by something like:

 \LabelsLocal

 \SetLabels
    (.5*.5) AAA\\
 \endSetLabels

 {
 \SetLabels
    (.5*.5) ZZZ\\
 \endSetLabels
   \AffixLabels{\FirstQuadrant}
 }

   \AffixLabels{\FirstQuadrant}

     There are however potential pitfalls.  Neither
labelfig.tex nor boxedeps.tex has been tested under extreme
conditions. Problems may occur if their procedures are
indiscriminately nested. For boxedeps.tex (not labelfig.tex)
there is a precise cause for worry, namely many of its
variables are "global", which means that TeX braces will not
provide the protection one might expect.

COMMAND SUMMARY FOR labelfig.tex

  Here [...] means optional (one or zero)
       [...]* means any number of such constructs

  \SetLabels
    [[<P>](<X><Sep><Y>) <label> \\]*
  \endSetLabels
  \ShowGrid  
  \AffixLabels{<the figure>}

   --- <P> is tack position, one of eleven or empty
              order irrelevant

                   \L\T      \T      \R\T

                   \L\E      \E      \R\E

                     \L               \R

                   \L\B      \B      \R\B

   --- (<X><Sep><Y>) insertion point;
  <Sep> is separator, = * by default;
  \SetXYSeparator{<Sep>} changes it.
   <X> and <Y> are real numbers

  --- <label> a label to attach 

  --- <the figure> the figure to label 

  \GlobalLabels (default)     
  \LocalLabels  setting for nested constructs.

 \Grids makes ALL grids appear; \HideGrid then makes just next disappear.
 \noGrids returns to default.  The commands are always global.

 \GridLineWidth{<dimension>} adjusts width of grid lines. Default is very
small, to give "hairline" effect. If your grid lines are missing try
setting \GridLineWidth{1pt}.

 \Edges#1 globally changes the edge size of all grids to the numerical 
value #1, which must be 1, 10, 100, or 1000.  The default is 1.

VERSION HISTORY.
 --- Jan 1993: \Edges#1 and [??] option after \SetLabels
 --- July 1992: \Grids, \noGrids, \HideGrid;
       Gridlines become hairlines; \GridLineWidth{<dimension>}.
 --- Oct 1991, Jan 1992: \SetXYSeparator{<Sep>},  \LabelsGlobal,
       \LabelsLocal.
 --- July 1991: first release

Address for bugs and other feedback:

        Raymond S\'eroul
        IREM and Lab. de Typographie Informatise
        Univ. Rene Descartes
        Strasbourg

    Tel 33-88-41-63-45
    Email:  A18645@FRCCSC21.BITNET

        Laurent Siebenmann
        Mathematique, Bat. 425,
        Univ de Paris-Sud,
        91405-Orsay,
        France

    Tel 33-1-6941-7949; 
    Email: lcs@topo.math.u-psud.fr

\title{%
  Determining Biholomorphic Type of a Manifold\\
  Using Combinatorial and Algebraic Structures %
}

\author{Sergiy A. Merenkov}{Merenkov, Sergiy A.}

\degree{Doctor of Philosophy}{Ph.D.}{August}{2003}

\majorprof{Alexandre Eremenko}

\newcommand\C{{\mathbb C}}
\newcommand\D{{\mathbb D}}
\newcommand\Hp{{\mathbb H}}
\newcommand\Z{{\mathbb Z}}
\newcommand\R{{\mathbb R}}
\newcommand\dee{\partial}
\newcommand\OC{{\overline{\mathbb C}}}

\newtheorem{fact}{Fact}
\newtheorem{remark}{Remark}
\newtheorem{lemma}{Lemma}

\begin{document}

%
%
%
%
\maketitle

\begin{dedication}
  To my teachers, Alex Eremenko at Purdue,\\
and Anatolij F. Grishin in Kharkov.
\end{dedication}

\begin{acknowledgments}
  I would like to express my thanks to my thesis advisor,
Alex Eremenko, who
provided me with very interesting problems, sharing his enlightening
ideas. It is my pleasure to thank David Drasin for his never-ending
encouragement and for his guidence through out the years. I would like to
thank Steve Bell for his mathematical discussions and his generous
support during a numerous semesters.
I thank Oded Schramm for giving his ideas concerning a counterexample
to R.~Nevanlinna's conjecture.
I also thank Burgess Davis for taking
part in the broadening  of my mathematical knowledge.
My general thanks to the faculty of the Department of Mathematics at
Purdue University for their help and support, especially to
L.~Avramov, D.~Catlin, A.~Gabrielov, and L.~Lempert.
I am thankful to all my
friends at Purdue for making my stay at Purdue as pleasant as it was.
\end{acknowledgments}

\newpage
\vspace*{2cm}
\begin{center}
This page deliberately left blank
\end{center}


\tableofcontents

\newpage
\vspace*{2cm}
\begin{center}
This page deliberately left blank
\end{center}


\listoffigures

\newpage
\vspace*{2cm}
\begin{center}
This page deliberately left blank
\end{center}





\begin{abstract}
We settle two problems of reconstructing a biholomorphic type of a
manifold. In the first problem we use graphs associated
to Riemann surfaces of a particular class. In the second one we use
the semigroup structure of analytic endomorphisms of domains in $\C^n$.

1. We give a new proof of a theorem due to P.~Doyle.
The problem is to determine a conformal type of a Riemann surface of
class $F_q$, using properties of the associated Speiser graph.
Sufficient criteria of type have been given since 1930's when
the class $F_q$ was introduced. Also there were necassary and sufficient
results which have theoretical value, but which are hard to apply.

P.~Doyle's theorem states that
a non-compact Riemann surface of class $F_q$ has a hyperbolic (parabolic)
type, if and only if its extended Speiser graph is hyperbolic (parabolic).
By a hyperbolic graph we mean a locally-finite infinite connected graph,
which admits a non-constant positive superharmonic function with respect
to the discrete Laplace operator. Otherwise a graph is parabolic.
The usefulness of this criterion stems from the possibility of applying
Rayleigh's short-cut method for graphs.

We apply Doyle's theorem to give a counterexample to a conjecture of
R.~Nevanlinna that relates the type to an excess of a Speiser graph.
More explicitely, the conjecture was that if the (upper) mean excess of a
surface of class $F_q$ is negative, then the surface is hyperbolic.
We provide an example of a parabolic surface of class $F_q$ with
negative mean excess.

2. If there is a biholomorphic or antibiholomorphic map
between two domains in $\C^n$, then it gives rise to an isomorphism
between the semigroups of analytic endomorphisms of these domains.

Suppose, conversely, that we are given two domains in $\C^n$ with
isomorphic semigroups of analytic endomorphisms. Are they
biholomorphically or antibiholomorphically
equivalent? This question was raised by L.~Rubel.
Similar questions were studied in the setting of topological spaces.

The case $n=1$ was investigated by A.~Eremenko, who showed that
if we require that the domains are bounded, then the answer to the
above question is positive. It was shown by A.~Hinkkanen that the boundedness
condition cannot be dropped.

We prove that two bounded domains in $\C^n$ with isomorphic
semigroups of analytic endomorphisms are biholomorphically or
antibiholomorphically equivalent.
Moreover, we generalize this by requiring only the existence of
an epimorphism between the semigroups.

\end{abstract}

%
%
%

\chapter{Introduction}

We study two problems, one of which deals with a class of Riemann surfaces 
represented by Speiser graphs, and the other one with bounded domains
in $\C^n$. Their settings and the methods we use to solve these problems
are different, but there is a unifying theme. Namely, in both cases
we determine a type,
conformal in the case of 
Riemann surfaces, or biholomorphic in the case of domains 
in $\C^n$, using an underlying combinatorial, respectively
algebraic structure. As an application to the first problem we give
an example showing that a conjecture of R.~Nevanlinna relating the
type of a surface to its excess is false.
A more detailed description of the problems follows. 

\section{P.~Doyle's Theorem}

A well-known theorem of Complex Analysis, the Uniformization Theorem,
says that every simply-connected Riemann surface is conformally equivalent 
to either the sphere, complex plane, or the unit disc. 
In the first case the surface is said to be of \emph{elliptic type}, in the 
second of \emph{parabolic type}, and in the third of \emph{hyperbolic type}.
When we come
up with a concrete Riemann surface, say by glueing together pieces
of the sphere along boundary parts, we would like to know how the
combinatorial pattern of glueing influences the type. One example of 
such a construction of Riemann surfaces is known in classical literature
as class $F_q$. These are the pairs $(X, f)$, where $X$ is a topological
manifold, and $f$ a continuous  open and discrete map from $X$ into the 
sphere $\OC$, so that $f$ is a covering map onto the sphere with
finitely many punctures. 

A surface of this class is uniquely represented
by a combinatorial object, called a Speiser graph, also known as 
a line complex, which is essentially the rule of pasting together
two complementary domains on the sphere, which share a Jordan curve as a 
common boundary. A Speiser graph is a homogeneous bipartite planar graph.
The components of its complement that are bounded by a finite number of edges
correspond to critical points, and those that are bounded by an infinite
number of edges correspond to assymptotic spots. 
Thus we come to the question of recovering the type from
properties of a Speiser graph. This problem has attracted a lot of
attention since the 1930's when the class $F_q$ was introduced. Many results 
relating properties of a graph to the type of the corresponding Riemann
surface have been obtained. Usually the criteria fall 
into one (and only one) of 
two categories: sharp but not useful, or useful but not sharp.

In 1984 Peter Doyle suggested a criterion of type which is sharp, and, at 
the same time,
seems to be useful (at least we were able to use it, unlike other 
known sufficient conditions, to provide a counterexample to 
R.~Nevanlinna's conjecture). The original proof due to Doyle, 
which is probabilistic
in nature, is very intuitive and enlightening, but might be hard to 
understand to non-specialists. 

Doyle's proof is based on the observation that the Brownian motion on 
a Riemann surface is transient if and only if there is a system of
currents out to infinity having finite dissipation rate. A {\it{system of 
currents out to infinity}} is a vector field, which is divergenceless 
outside of a sufficiently large compact set, and such that the total 
flux through
the boundary of this set is positive. The {\it{dissipation rate}}
of the flow is the integral of the square of the current density, i.e.
the square of the Hilbert-space norm of the vector field.
Similarly, the random walk on a graph is transient if and only if 
there is a system of currents through the edges of the graph out 
to infinity having finite dissipation rate. The dissipation rate
in this case is the sum of the squares of the currents through the edges.
Now, to prove the theorem one needs to show how a system of currents could be 
transfered from the surface to the associated graph and vice versa, 
without destroying the finiteness of the dissipation rate (see \cite{hR52}
for similar arguments).

We supply a new proof of Doyle's theorem. The methods we use are geometrical, 
and rely on the results due to M.~Kanai that assert the stability of
type under rough isometries, when the underlying spaces have bounded 
geometry. In accordance with this result, we construct a suitable conformal
metric on a given surface so that the surface equipped with this
metric is roughly isometric to the extended Speiser graph, introduced
by Doyle. An obvious choice for the metric would be the pullback of
the spherical metric, but unfortunately the surface equipped with this metric
is not roughly isometric to neither the Speiser graph, nor the extended
Speiser graph. In fact, no pullback metric can be suitable, since the
orders of critical points are in general unbounded.        

In Section~\ref{S:Back2}, we give the definition of a class of surfaces spread 
over the sphere, formulate the type problem, and provide background
information on graphs, Riemannian surfaces, and rough isometries. In
Section~\ref{S:Cl}, we give a definition of the class $F_q$ and examples.
Speiser graphs are introduced in Section~\ref{S:Spgr}. 
In Section~\ref{S:Doyle}, the extended
Speiser graph is introduced and the formulation of Doyle's theorem is given.
Section~\ref{S:ProofDoyle}, is devoted to the proof of Doyle's theorem.

\section{R.~Nevanlinna's Conjecture}

We give a counterexample to a conjecture of
R.~Nevanlinna that relates the type to the excess of a graph. 

For a Speiser graph $\Gamma$, R.~Nevanlinna introduces the
following characteristic. Let $V\Gamma$ denote the set of
vertices of the graph $\Gamma$. To each vertex $v\in V\Gamma$ we
assign the excess
$$
E(v)= 2-\sum_{f:\ v\in Vf} (1-1/k),
$$
where $f$ is a face with $2k$ edges, $k=1, 2, \dots, \infty$, and
$Vf$ is the set of vertices on its boundary. 
This notion is motivated via integral curvature, and thus reflects 
the geometric properties of the surface.

Nevanlinna also defines the mean excess of a Speiser graph
$\Gamma$. We fix a vertex $v\in V\Gamma$, and consider an
exhaustion of $\Gamma$ by a sequence of finite graphs
$\Gamma_{(i)}$, where $\Gamma_{(i)}$ is the ball of combinatorial
radius $i$, centered at $v$. By averaging $E$ over all the
vertices of $\Gamma_{(i)}$, and taking the limit, we obtain the
{\emph{mean excess}}, if the limit exists. We denote it by $E_m$.
If the limit does not exist, we consider {\it{upper}} or {\it
{lower excess}}, given by the upper, respectively lower, limit.
The upper mean excess of every infinite Speiser graph is $\leq 0$.

R.~Nevanlinna suggested a conjecture (\cite{rN70}, p. 312) that
a surface $(X, f)$ of the class $F_q$ is of a hyperbolic or a
parabolic type, according to whether the angle geometry of the surface is
``Lobachevskyan'' or ``Euclidean'', i.e. according to whether the mean
excess $E_m$ is negative or zero.

O.~Teichm\"uller gave an example of a surface of the hyperbolic type, 
whose mean excess is zero, thus contradicting a part of Nevanlinna's 
conjecture. We supply three examples contradicting the other part of the
conjecture, i.e. we construct parabolic surfaces with negative mean excess.
In the first example of a surface $(X, f)$, the function $f$ is analytic,
and in the second and third, $f$ does not have asymptotic values.
Thus we prove the following theorem.
\begin{theorem}\label{T:Counter}
There exists a parabolic surface $(X, f)\in F_3$ for which the
upper mean excess is negative.
\end{theorem}

In Section~\ref{S:Back3}, we recall definitions of the excess and the mean 
excess, illustrate these notions using integral curvature, and review
extremal length. In Sections~\ref{S:Ce1},~\ref{S:Ce2}, and~\ref{S:construct}, 
we provide the counterexamples. In Section~\ref{S:nonpos}, we 
construct an example of a simply connected,
complete, parabolic surface of nowhere positive
curvature, and such that its integral curvature in a disc around a 
fixed point is less than $-\epsilon$ times the area of the disc, for some
$\epsilon>0$ independent of the radius of the disc.

\section{Analytic Endomorphisms}

A classical theorem of L.~Bers says that every $\C$-algebra isomorphism 
$H(A)\rightarrow H(B)$ of algebras of holomorphic functions in domains $A$ 
and $B$ in the complex plane has either the form $f\mapsto f\circ\theta$, 
where 
$\theta:\ B\rightarrow A$ is a conformal isomorphism, or $f\mapsto 
\overline{f}\circ\theta$ with anticonformal $\theta$. In particular, the 
algebras $H(A)$ and $H(B)$ are isomorphic if and only if the domains $A$ 
and $B$ are conformally or anticonformally equivalent. 
H.~Iss'sa \cite{hI66} obtained a similar 
theorem for fields of meromorphic functions on Stein spaces. A good 
reference for these results is 
\cite{mH68}.

In 1990, L.~Rubel asked whether similar results hold for semigroups (under 
composition) $E(D)$ of holomorphic endomorphisms of a domain $D$. 
A question of recovering a topological space from the algebraic 
structure of its semigroup of continuous self-maps has been extensively
studied \cite{kM75}.

A.~Hinkkanen constructed 
examples \cite{aH92} which show that even non-homeomorphic domains in 
$\C$ can have isomorphic semigroups of endomorphisms. 
An elementary counterexample is a plane with 3 points removed 
and a plane with 4 points removed. They are obviously not biholomorphically
equivalent (they are not even homeomorphic for that matter), but if the
removed points are in general position, the corresponding semigroups 
consist of the unit and constant maps, and hence isomorphic.
The reason for such examples is 
that the semigroup of endomorphisms of a domain can be too small to 
characterize this domain. 

However, in 1993, A.~Eremenko \cite{aE93} proved that for two Riemann 
surfaces $D_1$, $D_2$, which admit bounded nonconstant holomorphic 
functions, and such that the semigroups of analytic endomorphisms $E(D_1)$ 
and $E(D_2)$ are isomorphic with an isomorphism 
$\varphi: \ E(D_1)\rightarrow E(D_2)$, there exists a conformal or 
anticonformal map $\psi:\ D_1\rightarrow D_2$ such that 
$\varphi f=\psi\circ f\circ\psi^{-1}$, for all $f\in E(D_1)$. 
We investigate the analogue of this result 
for the case of bounded domains in $\C^n$. The theorems of Bers 
and Iss'sa, mentioned above, do not extend to arbitrary domains in 
$\C^n$.

For a bounded domain $\Omega$ in $\C^n$ we denote by $E(\Omega)$ the 
semigroup of analytic endomorphisms of $\Omega$ under composition. We 
will write that a map is {\it{(anti-) biholomorphic}}, if it is 
biholomorphic 
or antibiholomorphic. We prove that 
if $\Omega_1,\ \Omega_2$ are bounded domains in $\C^n,\ \C^m$ 
respectively, and there exists 
$\varphi:\ E(\Omega_1)\rightarrow E(\Omega_2)$,  an isomorphism of 
semigroups, then $n=m$ and there exists an (anti-) 
biholomorphic map $\psi:\ \Omega_1\rightarrow\Omega_2$ such that 
\begin{equation}\label{E:C}
\varphi f=\psi\circ f\circ \psi^{-1}, \ \ \text{for all}\ f\in E(\Omega_1). 
\end{equation}

The existence of a homeomorphism $\psi$ follows from 
simple general considerations (Section~\ref{S:Top}). The hard part is proving 
that 
$\psi$ is (anti-) biholomorphic. In dimension 1 this is done by 
linearization of holomorphic germs of $f\in E(\Omega)$ near an attracting 
fixed point. In higher dimensions such linearization theory exists 
(\cite{vA88}, pp. 192--194), but it is too complicated (many germs with 
an attracting fixed point are non-linearizable, even formally). 
In Sections~\ref{S:Loc},~\ref{S:Ext}, we show how to localize the problem. 
In Sections~\ref{S:Sys},~\ref{S:Sim}  
we describe, using only the semigroup structure, a large enough class of 
linearizable germs. Linearization of these germs permits us to reduce the 
problem to a matrix functional equation, which is solved in 
Section~\ref{S:Sol}. In Section~\ref{S:Pr}, we complete the proof that
$\psi$ is (anti-) biholomorphic.

The above mentioned result can be slightly generalized, namely one may assume
that $\varphi$ is an epimorphism. In Section~\ref{S:Pr2}, we prove that
if $\varphi:\ E(\Omega_1)\rightarrow E(\Omega_2)$ is an epimorphism
between semigroups, where $\Omega_1,\ \Omega_2$ are bounded domains
in $\C^n,\ \C^m$ respectively, then $\varphi$ is an isomorphism.

\chapter{P.~Doyle's Theorem}

In this chapter we give an alternative proof of a theorem due to 
P.~Doyle~\cite{pD84a} on the type of a Riemann surface of class $F_q$. 

\section{Background and Preliminaries}\label{S:Back2}

\subsection{Uniformization Theorem}

A \emph{Riemann surface} is a 1-dimensional complex manifold, 
or, in other words,
it is a 2-real-dimensional manifold endowed with a maximal atlas in
which all transition maps are conformal. It is \emph{simply-connected} if
the fundamental group is trivial.

The following well-known fact is called the Uniformization Theorem 
\cite{lA73}.
\begin{theorem}  
For every simply-connected
Riemann surface $X$ there exists a conformal map $\varphi:\ X_0\to
X$, where $X_0$ is one of the three model surfaces:
\begin{enumerate}
\item 
the open unit disc $\D_1$;
\item
the complex plane $\C$;
\item
the extended complex plane $\OC$. 
\end{enumerate}
\end{theorem}
The map $\varphi$ is called the {\emph{uniformizing
map}}. The Uniformization Theorem has a number of applications, the main 
of which is that on every Riemann surface there exists a conformal
metric of constant Gaussian curvature -1, 0, or 1. 
\begin{definition}
A simply-connected Riemann surface $X$ is said to have a
hyperbolic, parabolic, or elliptic
type, according to whether it is conformally equivalent to $\D_1$,
$\C$, or $\OC$ respectively.
\end{definition}
Sometimes we simply say that $X$ is hyperbolic, parabolic, or elliptic. 
Also, we refer to the type of a simply-connected Riemann surface 
as a conformal type.

\subsection{Surfaces Spread over the Sphere}

We are interested in the application of the Uniformization Theorem
to the following construction. A map between two topological spaces is
called \emph{open}, if the image of every open set is open. It is called
\emph{discrete}, if the preimage of every point is discrete, i.e.
every point of the preimage has a neighborhood that does not contain any other
points of the preimage.
\begin{definition}
A surface spread over the
sphere is a pair $(X, f)$, where $X$ is a topological surface
and $f:\ X\to\OC$ a continuous, open and discrete map. 
\end{definition}
The map $f$ is called a {\emph{projection}}. 
Two such surfaces $(X_1, f_1),\
(X_2, f_2)$ are {\emph{equivalent}}, if there exists a
homeomorphism $\phi:\ X_1\to X_2$, such that $f_1= f_2\circ\phi$.
A theorem of Sto\"{\i}low \cite{sS56} implies that for every continuous
open and discrete 
map $f$ from a topological surface (i.e. a 2-real-manifold) 
to the Riemann sphere there exists a homeomorphism $\phi$ of $X$ onto
a Riemann surface $Y$, so that the map $f\circ\phi^{-1}$ is meromorphic. 
The Riemann surface $Y$ is unique up to conformal equivalence. This tells us
that there exists a unique conformal structure on $X$ (i.e.
$X$ becomes a Riemann surface), which makes $f$ into a meromorphic
function. Near each point $x\in X$ 
the function $f$ is conformally equivalent to a
map $z\mapsto z^k$, with $k$ depending on $x$. The number $k=k(x)$
is called the {\emph{local degree}} of $f$ at $x$. If $k\neq 1$,
$x$ is called a {\emph{critical point}} and $f(x)$ a
{\emph{critical value}}. The set of critical points is a discrete
subset of $X$.

The surface $X$ can be endowed with a metric that is the 
$f$-pullback of the
spherical metric $2|dw|/(1+|w|^2)$. The pullback metric is
singular, i.e. it is degenerate on a discrete set in $X$. The
surface $X$, endowed with the pullback metric, is a particular
case of spherical polyhedral surfaces \cite{mB00a}, \cite{mB02}.

\subsection{Type Problem}

If $X$ is simply-connected, what is the type
of the Riemann surface obtained as in the previous section,
if $(X, f)$ is a surface spread over the sphere? More precisely,
how does the conformal type depend on the properties of the function
$f$ that are invariant under homeomorphic changes of the independent variable? 
This is the formulation of the type problem.

By uniqueness of the conformal structure, 
equivalent surfaces have the same type. We notice
that it is easy to single out the elliptic type as consisting of
compact Riemann surfaces. So we are down to the choice between
hyperbolic and parabolic types.

The dependence of type on curvature properties has been studied in
\cite{lA73}, 
\cite{mB00}, \cite{yR93}. We study the type 
problem for surfaces of so called class $F_q$ in Section~\ref{S:ClF}. To
surfaces of this class one can naturally associate a planar graph, 
Section~\ref{S:Spgr}, called a Speiser graph. We are interested in 
the dependence of type of properties of this graph.

\subsection{Graphs}

By a \emph{graph} $G$ we mean a pair $(V, E)$, where $V$ is 
an at most countable set, whose elements are called \emph{vertices},
and $E$ a set of pairs of elemets from $V$. Elements of $E$ are
called \emph{edges}. We say that $e\in E$ \emph{connects} 
$v_1, v_2\in V$, or that $e$ is an edge between $v_1$ and $v_2$, if
$e=(v_1, v_2)$. Multiple edges  between two vertices are allowed, 
but loops, i.e. edges of the form $(v, v)$ are not.
  
Given a connected graph $G$, we denote by $VG$, $EG$ the sets of
its vertices and edges respectively. If two vertices $v_1, v_2$ of $G$ are
connected by an edge, we write $v_1\sim v_2$. We denote by
$\deg_vG$, the number of edges of $G$ emanating from $v$. A graph
$G$ is said to have a {\emph{bounded degree}}, if $\sup\{\deg_vG:
\ v\in VG\}<\infty$. If $G'$ is a connected subgraph of $G$, the
{\emph{boundary}} of $G'$ is the set of vertices $v\in VG'$, such
that $\deg_vG'<\deg_vG$. A {\emph{path}} in $G$ is a connected
subgraph, which has degree 2 at all of its vertices with at most
two exceptions, where it has degree 1. A connected graph $G$ is a
metric space with a combinatorial distance on it, i.e. the
distance between two vertices is the number of edges of a shortest
path connecting them. If $G$ is a connected graph embedded in a
topological surface $X$, the connected components of $X\setminus
G$ are called {\emph{faces}} of $G$; the set of faces of $G$ is
denoted by $FG$. For a graph $G$ embedded in the plane, we denote
by $G^*$ its dual.

If a graph $G$ is locally-finite, a linear operator $\Delta$,
acting on functions $u$ on $VG$, is defined by
$$
\Delta u(v)=1/\deg_vG\sum_{v'\sim v}u(v')-u(v), \ \ v\in VG.
$$
It is well-known that the operator $\Delta$ enjoys many properties
that the Laplace operator possesses \cite{eD69}. A locally-finite
infinite connected graph $G$ is called {\emph{hyperbolic}}, if
there exists a positive non-constant superharmonic function on
$VG$. Otherwise it is called {\emph{parabolic}}. The hyperbolicity
(parabolicity) for a locally-finite infinite graph is equivalent
to the transience (recurrence) of the simple random walk on it.

\subsection{Riemannian Surfaces}

A {\emph{conformal metric}} $ds$ on a Riemann surface $X$ is a
metric whose length element is given in local coordinates by
$\rho(z)|dz|$, where $\rho$ is a positive smooth function. Often
one considers a more general conformal metric, by allowing $\rho$
to vanish on a discrete set. For example, spherical polyhedral
surfaces mentioned above carry such a metric. For our purposes
conformal metrics with everywhere positive $\rho$ will be
sufficient.

The \emph{Gaussian curvature} of a conformal metric $\rho(z)|dz|$
is given by $-\rho(z)^{-2}\delta\log\rho(z)$. It is isometry 
invariant.

We denote by $Y=(X, ds)$ a pair, where $X$ is a Riemann surface
and $ds$ is a conformal metric on $X$. We call such a $Y$ a
{\emph{Riemannian surface}}. A Riemannian surface $Y$, not
necessarily simply-connected, is said to be {\emph{hyperbolic}},
if there exists a positive non-constant superharmonic function on
it. Otherwise it is called {\emph{parabolic}}. Since the metric is
conformal, a superharmonic function on $Y$ is the same as a
superharmonic function on $X$. Therefore, a simply-connected
Riemann surface $X$ is conformally equivalent to $\D_1$ ($\C$), if
and only if $Y=(X, ds)$ is hyperbolic (parabolic) as a Riemannian
surface with an arbitrary conformal metric $ds$ on it. Moreover,
the following fact holds.
\begin{fact}\label{F:fact}
If $A$ is an arbitrary discrete subset of an open simply-connected
Riemann surface $X$, and $ds$ is a conformal metric on $X\setminus
A$, then $Y=(X\setminus A, ds)$ is hyperbolic, if and only if $X$
is conformally equivalent to $\D_1$.
\end{fact}
This is because every positive superharmonic function on
$X\setminus A$ extends to a superharmonic function on $X$
\cite{wH76}.

A Riemannian surface is \emph{complete} if it is complete as
a metric space. A \emph{radius of injectivity} of a Riemannian
surface $Y$ is the infimum over all points $x$ of $Y$ of the 
supremum over all non-negative $r$ such that for all $t\leq r$
the ball centered at $x$ of radius $t$ is homeomorphic to a 
Euclidean ball.

We say that a Riemannian surface $Y$ satisfies the
{\emph{geometric uniformness condition}}, if
\begin{equation}\label{E:Guc}
\begin{split} Y {\text{ is complete, the Gaussian curvature is
bounded}}\\
{\text{from below, and the radius of injectivity is positive.}}
\end{split}
\end{equation}

\subsection{Rough Isometry}

Let $(X_1, d_1)$ and $(X_2, d_2)$ be two metric spaces.
\begin{definition}\label{D:Ris}
A map $\Phi:\ X_1\to X_2$, not necessarily continuous, is called a
rough isometry, if the following two conditions are satisfied:
\begin{enumerate}
    \item for some $\epsilon>0$, the $\epsilon$-neighborhood
    of the image of $\Phi$ in $X_2$ covers $X_2$;
    \item there are constants $C_1\geq 1, \ C_2\geq 0$, such that
    for all $x_1,\ x_2\in X_1$,
$$
C_1^{-1}d_1(x_1, x_2) - C_2\leq d_2(\Phi(x_1), \Phi(x_2))\leq
C_1d_1(x_1, x_2) + C_2.
$$
\end{enumerate}
\end{definition}

A metric space $(X_1, d_1)$ is said to be {\it{roughly isometric}}
to a metric space $(X_2, d_2)$, if there exists a rough isometry
from $X_1$ into $X_2$. This is an equivalence relation. The notion
of rough isometry was introduced by M.~Kanai \cite{mK85} and
M.~Gromov \cite{mG81}.

An immediate consequence of Kanai's results \cite{mK85},
\cite{mK86} is the following theorem.
\begin{theorem}\label{T:Mgt}
If a non-compact Riemannian surface $Y$ satisfying~(\ref{E:Guc})
is roughly isometric to a connected locally-finite graph $G$ of
bounded degree, then $Y$ is hyperbolic, if and only if $G$ is
hyperbolic.
\end{theorem}
In fact, Kanai proves that a Riemannian surface $Y=(X, ds)$ is
hyperbolic, if and only if an $\epsilon$-net in $Y$ is hyperbolic.
An {\emph{$\epsilon$-net}} in $Y$ is a maximal
$\epsilon$-separated set $Q$ in $D$ with a structure of a graph,
so that vertices are points of $Q$; two vertices $v_1, v_2$ are
connected by an edge, if and only if $d(v_1, v_2)\leq 2\epsilon$.
The graph $Q$ has a bounded degree and is roughly isometric to
$Y$. A graph $G$, roughly isometric to $Y$ is, by transitivity,
roughly isometric to $Q$. Since both graphs have bounded degree,
they are \cite{pS94} simultaneously hyperbolic or parabolic.

\section{Class $F_q$}\label{S:ClF}\label{S:Cl}

\subsection{Definition}

We study the type problem for a particular, but rather broad,
subclass of surfaces spread over the sphere, the so called class
$F_q$. For a surface of this class we investigate the dependance
of type on the properties of the associated Speiser graph (see~\ref{S:Spgr}). 
In this respect see \cite{pD84}, \cite{aG70}, \cite{rN70}, \cite{aS30},
\cite{oT38}, \cite{lT47}, \cite{lV50}, \cite{hW68}.

Let $\{a_1, \dots, a_q\}$ be distinct points in $\OC$.
\begin{definition}\label{D:Fq} A surface $(X, f)$, where $X$ is 
open and simply-connected,
belongs to class $F_q=F(a_1, \dots, a_q)$, if
$$
f:\ X\setminus\{f^{-1}(a_i),\ i=1, \dots, q\}\to
\OC\setminus\{a_1,\dots, a_q\}
$$
is a covering map.
\end{definition}

Analytically surfaces of class $F_q$ can be characterized as those for
which the function $f$ has only finitely many critical and asymptotic values.
An \emph{asymptotic spot} is an open arc contained in $X$ that escapes
from every compact subset of $X$, and such that the limit of $f$ 
along this arc exists.
An \emph{asymptotic value} is the limit of $f$ at an asymptotic spot.

For each $i$, let $(V_i, \psi_i)$ be a coordinate neighborhood of
$a_i$, centered at zero, so that $\psi_i(V_i)=\D_1$, and $V_i\cap
V_j=\emptyset,\ i\neq j$. The restriction of $f$ to a connected
component $U$ of $f^{-1}(V_i\setminus{a_i})$ is a covering map.
Therefore \cite{oF91} this map is conformally equivalent to either
$\D_1^*\to\D_1^*,\ z\mapsto z^k$, or $\Hp\to\D_1^*,\
z\mapsto \exp(z)$, where $\D_1^*$ denote the punctured open unit
disc, and $\Hp$ an open left half-plane. In particular, $U$ does
not contain any critical points of $f$.

\subsection{Examples}

The following are examples of surfaces of class $F_q$.
\begin{enumerate}
\item
$(\C, \sin)\in F_3(-1, 1, \infty)$.
\item
$(\D_1, \lambda)\in F_3(0, 1, \infty)$,
where $\lambda$ is a modular function.
\end{enumerate}

\section{Speiser Graphs}\label{S:Spgr}

\subsection{Definition}

We fix a Jordan curve $L$, containing the points $a_1,\dots, a_q$.
The curve $L$ is usually called a {\emph{base curve}}. It
decomposes the sphere into two simply-connected regions $H_1,\
H_2$, called {\emph{half-sheets}}. We assume that the indices of
$a_i$'s are cyclically ordered modulo $q$, and the curve $L$ is
oriented so that the region $H_1$ is to the left. We denote by
$L_i$ the arc on $L$ between $a_i$ and $a_{i+1}$. Let us fix
points $p_1$ in $H_1$ and $p_2$ in $H_2$, and choose $q$ disjoint
Jordan arcs $\gamma_1,\dots, \gamma_q$ in ${\OC}$, such that each
arc $\gamma_i$ has $p_1$ and $p_2$ as its endpoints, and has a
unique point of intersection with $L$, which is on $L_i$. Let
$\Gamma'$ denote the graph embedded in $\OC$, whose vertices are
$p_1,\ p_2$, and edges $\gamma_i,\ i=1,\dots, q$, and let $\Gamma$
be the $f$-pullback of the graph $\Gamma'$. We identify $\Gamma$
with its image in $\R^2$ under a sense-preserving homeomorphism of
$X$ onto $\R^2$. Clearly it does not depend on the choice of the
points $p_1,\ p_2$, and the curves $\gamma_i,\ i=1, \dots, q$. 
The graph $\Gamma$ has the following properties:
1. $\Gamma$ is infinite, connected,
2. $\Gamma$ is homogeneous of degree $q$, and
3. $\Gamma$ is bipartite.

A graph, properly embedded in the plane and satisfying properties
1, 2, and 3, is called a {\emph{Speiser graph}}, also known as a
{\emph{line complex}}. The vertices of a Speiser graph $\Gamma$
are traditionally denoted by $\times$ and $\circ$. Each face of
$\Gamma$, i.e. a connected component of $\R^2\setminus\Gamma$, has
either $2k$ edges $k=1,2,\dots$, in which case it is called an
{\emph{algebraic elementary region}}, or infinitely many edges,
called a {\emph{logarithmic elementary region}}. Two Speiser
graphs $\Gamma_1, \ \Gamma_2$ are said to be {\emph{equivalent}},
if there is a sense-preserving homeomorphism of the plane, which
takes $\Gamma_1$ to $\Gamma_2$. Below we refer to an equivalence
class as a Speiser graph.

\subsection{Examples}
\begin{figure}[h]
\centerline{\scalebox{.8}{\includegraphics{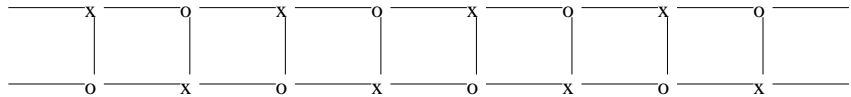}}}
\caption{Speiser graph of sine}
\label{fig:sine0}
\end{figure}
\begin{figure}[h]
\centerline{\scalebox{.85}{\includegraphics{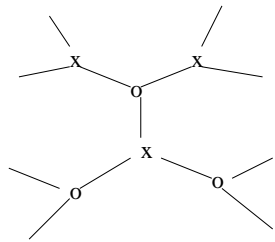}}}
\caption{Speiser graph of $\lambda$}
\label{fig:sine1}
\end{figure}

\subsection{Reconstructing a surface from a Speiser Graph}

The above construction of a Speiser graph of a surface $(X, f)\in F_q$
is reversible. Suppose that the faces of a
Speiser graph $\Gamma$ are labelled by $a_1, \dots, a_q$, so that
when going counterclockwise around a vertex $\times$, the indices
are encountered in their cyclic order, and around $\circ$ in the
opposite cyclic order. A labelling of faces induces the one of
edges: we assign a label $i$ to an edge, if it is the common
boundary for faces labelled $a_i$ and $a_{i+1}$. We fix a base
curve $L$ in $\OC$ passing through $a_1, \dots, a_q$ in the order
of increasing indices, and denote by $H_1$ and $H_2$ the
half-sheets, so that $H_1$ is to the left. Then one constructs a
surface $(X, f)\in F(a_1, \dots, a_q)$ in the following way. Let
$\Gamma^*$ be the cell decomposition of $\R^2$, dual to $\Gamma$.
Each 2-dimensional cell has $q$
1-dimensional cells on its boundary.
The 2-dimensional cells are labelled by $\times$ and $\circ$, and the
0-dimensional cells by $a_1, \dots, a_q$. We
map each 2-dimensional cell of $\Gamma^*$ labelled by $\times$ to $H_1$, and
each 2-dimensional cell labelled by $\circ$ to $H_2$, so that the maps agree on
common boundary 1-dimensional cells and a 0-dimensional cell labelled by $a_i$ 
is mapped to $a_i$.
Thus we obtain a continuous, open and discrete map $f:\
\R^2\to\OC$, such that $\R^2\setminus\{f^{-1}(a_i), i=1, \dots,
q\}\to \OC\setminus\{a_1, \dots, a_q\}$ is a covering map. So,
$(\R^2, f)\in F(a_1, \dots, a_q)$, and its Speiser graph is clearly
$\Gamma$.

It is natural to ask whether $(X, f)\in F_q$ is hyperbolic if and
only if its Speiser graph $\Gamma$ is hyperbolic. The intuition
behind this question is in viewing the simple random walk on
$\Gamma$ as a discrete approximation of the Brownian motion on the
Riemann surface \cite{sK49}, \cite{sK53}. Unfortunately, as we
show in Appendix~\ref{S:Noten}, the hyperbolicity of a surface of
the class $F_q$ is not equivalent to the hyperbolicity of its
Speiser graph.

\section{P.~Doyle's Theorem}\label{S:Doyle}

\subsection{Extended Speiser Graph}

P. Doyle \cite{pD84} suggested to use
an extended Speiser graph to study the type problem.

Let $\Z$ denote the set of integers, and $\Z_+$ the set of
non-negative integers.

A {\emph{half-plane lattice}} $\Lambda$ is the graph embedded in
$\R^2$, whose vertices form the set $\Z\times\Z_+$. Two vertices
$(x', y'),\ (x'', y'')$ are connected by an edge, if and
only if $(x''-x', x''-x')=(\pm 1, 0)$ or $(0, \pm 1)$. The
boundary of the half-plane lattice is the infinite connected
subgraph, whose set of vertices is $\Z\times\{0\}$. There is an
action of $\Z$ on $\Lambda$ by horizontal shifts. A
{\emph{half-cylinder lattice}} $\Lambda_n$ is $\Lambda/n\Z$. The
boundary of $\Lambda_n$ is the induced boundary from $\Lambda$.

Let $n\geq1$ be given. If we replace each face of a Speiser graph
$\Gamma$ with $2k$ edges, $k\geq n$, by the half-cylinder lattice
$\Lambda_{2k}$, and each face with infinitely many edges by the
half-plane lattice $\Lambda$, identifying the boundaries of the
faces with the boundaries of the corresponding lattices along the
edges, we obtain the {\emph{extended Speiser graph}} $\Gamma_n$.
The graph $\Gamma_n$ is an infinite connected graph, embedded in
$\R^2$, and containing $\Gamma$ as a subgraph. It has a bounded
degree, and all faces of $\Gamma_n$ have no more than $\max\{2(n-1), 4\}$
edges.

\subsection{Statement of the Theorem}

\begin{theorem}\label{T:Ut}
For every  $n\geq 1$, a surface $(X, f)\in F_q=F(a_1, \dots, a_q)$
has a hyperbolic (parabolic) type, if and only if $\Gamma_{n}$ is
hyperbolic (parabolic), where $\Gamma$ is the Speiser graph of
$(X, f)$.
\end{theorem}

Theorem~\ref{T:Ut} is a slight generalization of 
the theorem by P. Doyle \cite{pD84}. The latter
states that $(X, f)\in F_q$ is hyperbolic, if and only if the
McKean-Sullivan random walk on its Speiser graph $\Gamma$ is
transient. In plain terms, the McKean-Sullivan random walk on
$\Gamma$ comes from the simple random walk on $\Gamma_1$, when we
observe it only as it hits $\Gamma$. Doyle's arguments are
probabilistic and electrical, whereas we employ geometric methods.
We use the results of M. Kanai \cite{mK85}, \cite{mK86} to prove
Theorem~\ref{T:Ut}.

\section{Proof of P.~Doyle's Theorem}\label{T:T2}\label{S:ProofDoyle}

According to Fact~\ref{F:fact} and Theorem~\ref{T:Mgt}, to prove
Theorem~\ref{T:Ut}, we need to find a conformal metric $ds$ on
$X\setminus A$, where $A$ is a discrete subset of $X$, such that
$Y=(X\setminus A, ds)$ satisfies~(\ref{E:Guc}), and is roughly
isometric to $\Gamma_n$.

\subsection{Conformal Metric}

For each $i=1, \dots, q$, let $(V_i, \psi_i)$ be a local
coordinate neighborhood of $a_i$, centered at zero, so that
$\psi_i(V_i)=\D_1$, and $\overline{V}_i\cap\overline{V}_j=\emptyset,
\ i\neq j$.
Consider an open covering $\{V_i:\ i=0, \dots, q\}$ of $\OC$,
where $V_0=\OC\setminus\{\overline{\psi_i^{-1}(\D_{1/2})}:\ i=1,
\dots, q\}$; let $\psi_0$ be a conformal map of $V_0$ onto a
domain in $\C$. Further, let $\{g_i,\ i=0, \dots, q\}$ be a
partition of unity on $\OC$, subordinate to the covering. This
partition of unity pulls back to a partition of unity on $X$ as
follows. Let $W$ be a connected component of $f^{-1}(V_i)$. We
define $g_W=g_i\circ f$, a function on $W$, that we extend to a
smooth function on $X$, by letting it to be 0 outside $W$. It is
clear that $\{g_W\}$, a family of functions indexed by connected
components of $f^{-1}(V_i),\ i= 0, 1, \dots, q$, forms a partition
of unity on $X$. Every component $W$ contains at most one singular
point. If $(W, f)$ is a $k$-sheeted covering of $V_i$, $k=1, 2,
\dots, \infty$, we denote $W$ by $W_k$. The connected component of
$f^{-1}(V_0)$ is denoted by $W_0$.

We choose $A$ to be the set of all critical points $x\in X$, so
that the local degree $k=k(x)$ of $f$ at $x$ is at least $n$. Now we
define a conformal metric $ds$ on $X\setminus A$:
$$
\begin{aligned}
\rho(z)|dz|&=g_{W_0}f^*(\psi_0^*(|dw|))+
\sum_{i>0}\bigg\{\sum_{W_k:\ k\geq n}
g_{W_k}f^*(\psi_i^*(|dw|/|w|))\\
&+\sum_{W_k:\ k< n}
g_{W_k}f^*(\psi_i^*(|dw|/|w|^{(k-1)/k}))\bigg\}.
\end{aligned}
$$
The function $\rho$ smoothly extends to a neighborhood of every
critical point that does not belong to $A$.

We need to show that $Y=(X\setminus A, ds)$ satisfies
$(\ref{E:Guc})$, and is roughly isometric to $\Gamma_{n}$.

\subsection{Geometric Uniformness Condition}

Every curve going out to a point in $A$ has infinite length in the
metric $ds$, thus $Y$ is complete. The Gaussian curvature is
bounded. Indeed, suppose that it is not. Then there exists a
sequence of points in $X\setminus A$, on which the Gaussian curvature tends
to $\infty$. We project this sequence to $\OC$. The projected
sequence has either finitely many points, or accumulates to a
point in $\OC$. Since in a neighborhood of each point in $\OC$
there are at most $n+1$ choices for the metric, each with a bounded
curvature, we get a contradiction.

The radius of injectivity is positive. Assume the contrary, i.e.
there exists a sequence of points $\{x_n\}$ in $X\setminus A$,
such that if $r_n$ is the radius of injectivity at $x_n$, then
$r_n\to 0$. To derive a contradiction, we follow the same argument
as in the proof that the Gaussian curvature is bounded. The most
interesting case is when the projected sequence accumulates at a
point $a_i$. Every component $W_k'\subset W_k,\ k<\infty$, of the
preimage $\psi_i^{-1}(\D_{1/2})$ is isometric to
$\D_{(1/2)^{1/k}}$ with the length element $k|dz|/|z|$, if $k\geq
n$, or $k|dz|$, when $k<n$. A connected component
$W_{\infty}'\subset W_{\infty}$ of $\psi_i^{-1}(\D_{1/2})$ is
isometric to $\{z, \Re z<1/2\}$ with the metric $|dz|$. In any
case we obtain a contradiction.

\subsection{Rough Isometry}

It remains to show the rough isometry. The Speiser graph $\Gamma$
of $(X, f)$ is the preimage of $\Gamma'$ under $f$, where
$\Gamma'$ is embedded in $\OC\setminus \{a_i,\ i=1, \dots, q\}$.
The graph $\Gamma'$ is finite, and otherwise satisfies all the
properties that a Speiser graph does. Therefore we can form the
extended graph $\Gamma_1'$. Since each face of $\Gamma'$ contains
a unique $a_i$, we can assume that the extended graph $\Gamma_1'$
is embedded in $\OC\setminus\{a_1, \dots, a_q\}$ in such a way,
that with respect to the local coordinate $(V_i, \psi_i)$, the
edges of the lattice of $\Gamma_1'$ are Euclidean semicircles and
orthogonal to them family of straight segments, which have length
1 in the metric $\psi_i^*(|dw|/|w|)$.

Let $J:\ \Gamma_{n}\to X\setminus A$ be the embedding, whose image
is contained in the pullback of $\Gamma_1'$. For this embedding
properties 1, 2, and 3 of Lemma~\ref{L:Ri} below are readily
verified, using the fact that there is a finite number of choices
for the metric in $f(X\setminus A)$. Theorem~\ref{T:Ut} follows.

For a Riemannian surface $Y$ we denote by $d_Y(x_1, x_2)$ the
distance between $x_1, x_2 \in Y$, and by $l_Y(C)$, the length of
a curve $C\subset Y$. Similarly, for a graph $G$, we denote the
combinatorial distance between $v_1, v_2 \in VG$ by $d_{G}(v_1,
v_2)$, and the combinatorial length of a path $C$ by $l_{G}(C)$. A
curve in $Y$ joining points $x_1$ and $x_2$ is denoted by $C_{x_1,
x_2}$.

\begin{lemma}\label{L:Ri}
Suppose that for a connected graph $G$ properly embedded in a
complete Riemannian surface $Y$  the following conditions are
satisfied:
\begin{enumerate}
\item
there exists a constant $\epsilon$, such that for every point
$x\in Y$,
$$
\inf\{d(x, v):\ v\in VG\}<\epsilon,
$$
\item
there exist constants $B_1, B_2,\ 0<B_1<B_2$, such that for every
edge $e\in EG$,
$$
B_1\leq l_{Y}(e)\leq B_2,\ \ {\text{and}} 
$$
\item
there exists a constant $B_3>0$, such that for every face $f\in
FG$, and every two points $x_1, x_2$ on the boundary $\dee f$ of
$f$,
$$
\inf\{l_{Y}(C_{x_1, x_2}):\ C_{x_1, x_2}\subset f\} \geq B_3
\inf\{l_{Y}(C_{x_1, x_2}):\ C_{x_1, x_2}\subset \dee f \}.
$$
\end{enumerate}
Then  the graph $G$ is roughly isometric to $Y$.
\end{lemma}

\emph{Proof.} Let $J:\ G\to Y$ be the embedding map. In view of
condition {\it{1}}, the first property of rough isometry for $J$ is
satisfied, so it remains to prove the second property.

Let $v_1, v_2 \in VG\subset Y$, and $C$ be a path in $G$, joining
these two points, and having the minimal combinatorial length.
Then, by condition~{\it{2}},
\begin{equation}\label{E:Os}
d_{G}(v_1, v_2)=l_{G}(C)\geq (1/B_2) l_{Y}(C)\geq
(1/B_2)d_{Y}(v_1, v_2).
\end{equation}

Conversely, let $v_1, v_2\in VG\subset Y$, and $C$ be a curve in
$Y$, joining $v_1$ and $v_2$. Let $f$ be a face of $G$, such that
$C\cap f\neq \emptyset$, and $C_f$ be a curve which is a connected
component of $C\cap f$. If $x_1, x_2$ are endpoints of $C_f,\ x_1,
x_2\in \dee f$, then, by condition~{\it{3}},
$$
l_{Y}(C_f)\geq B_3 \inf\{ l_{Y}(C_{x_1, x_2}):\ C_{x_1,
x_2}\subset \dee f \}.
$$
Since this holds for every face $f$ and every component $C_f$, we
conclude that there exists a path $C'$ in $G$, joining $v_1$ and
$v_2$, such that
$$
l_{Y}(C)\geq B_3l_{Y}(C')\geq B_1B_3 d_{G}(v_1, v_2),
$$
where the last inequality holds by condition~{\it{2}}. Taking the infimum
with respect to curves $C$ joining $v_1$ and $v_2$, we obtain that
\begin{equation}\label{E:S3}
d_{Y}(v_1, v_2)\geq B_1 B_3 d_{G}(v_1, v_2).
\end{equation}
Combining inequalities (\ref{E:Os}), (\ref{E:S3}),
$$
B_1B_3 d_{G}(v_1, v_2)\leq d_{Y}(v_1, v_2)\leq B_2 d_{G}(v_1,
v_2),
$$
we conclude that $J$ is a rough isometry. The lemma is proved. $\Box$

An immediate corollary of Lemma~\ref{L:Ri} is the fact that a
surface of the class $F_q$, endowed with the pullback of a
spherical metric, is roughly isometric to the dual of its Speiser
graph. However we could not use this rough isometry in studying
the type problem due to the presence of vertices of infinite
degree on the dual of a Speiser graph. Also, even if we assume that 
there are no asymptotic values, i.e. there are no vertices of infinite
degree on a dual of the Speiser graph, the degrees of the vertices
of the dual can be unbounded, and we cannot conclude that the type
of a surface agrees with the type of the dual graph. 

The only case when we can use the dual graph to determine the type
of a surface is when the degrees of the vertices of this graph are bounded. 
This is first of all too restrictive, and second, if this happens,
we can use the Speiser graph itself for this purpose, i.e. we
do not need to consider the extended graph. 

In the next chapter we apply Doyle's theorem to show that R.~Nevanlinna's
conjecture is false. We could not use any other known criteria of type
to show the parabolicity of the surface constructed below.

\newpage
\vspace*{2cm}
\begin{center}
This page deliberately left blank
\end{center}

\chapter{R.~Nevanlinna's Conjecture}

In this chapter we provide three examples of a parabolic surface with
negative mean excess, contradicting Nevanlinna's conjecture. We also 
give an example of a parabolic surface with ``a lot of negative curvature'', 
Section~\ref{S:nonpos}.
The third example and the surface in Section~\ref{S:nonpos}, are
due to O.~Schramm and I.~Benjamini.

\section{Background and Preliminaries}\label{S:Back3}

\subsection{Excess}

For a Speiser graph $\Gamma$, R.~Nevanlinna introduces the
following characteristic. To each vertex $v\in V\Gamma$ we assign
the number
$$
E(v)= 2-\sum_{f:\ v\in Vf} (1-1/k),
$$
where $f$ is a face with $2k$ edges, $k=1, 2, \dots, \infty$, and
$Vf$ is the set of vertices on its boundary. The function $E:\
V\Gamma\to \R$, $v\mapsto E(v)$ is called the {\emph{excess}} of
$\Gamma$. 

\subsection{Integral Curvature}

The motivation for the definition of excess uses a notion 
of integral curvature.  The {\emph{integral curvature}} $\omega$ on $X$ is
a signed Borel measure, so that for each Borel subset $B\subset
X$, $\omega(B)$ is the area of $B$ with respect to the pullback
metric minus $2\pi\sum(k-1)$, where the sum is over all critical
points $x\in B$, and $k$ is the local degree of $f$ at $x$.

Each vertex of $\Gamma$ represents a hemisphere, and each face of
$\Gamma$ with $2k$ edges, $k=2, 3, \dots$, represents a critical
point, where $k$ is the local degree of $f$ at this point. Therefore,
each vertex of $\Gamma$ has positive integral curvature $2\pi$,
and each face with $2k$ edges has negative integral curvature
$-2\pi(k-1)$. We spread the negative curvature evenly to all the
vertices of the face. A face with infinitely many edges
contributes $-2\pi$ to each vertex on its boundary. The curvature
mass obtained by every $v\in V\Gamma$ is exactly $2\pi E(v)$.

\subsection{Mean Excess}

Nevanlinna also defines the mean excess of a Speiser graph
$\Gamma$. We fix a vertex $v\in V\Gamma$, and consider an
exhaustion of $\Gamma$ by a sequence of finite graphs
$\Gamma_{(i)}$, where $\Gamma_{(i)}$ is the ball of combinatorial
radius $i$, centered at $v$. By averaging $E$ over all the
vertices of $\Gamma_{(i)}$, and taking the limit, we obtain the
{\emph{mean excess}}, if the limit exists. We denote it by $E_m$.
If the limit does not exist, we consider the {\it{upper}} or {\it
{lower excess}}, given by the upper, respectively lower, limit.
The upper mean excess of every infinite Speiser graph is nonpositive (see 
Appendix~\ref{Upmeanex}).

\subsection{Extremal length}\label{S:Er}

In this section we give the definition of the extremal length of a
family of paths and derive one of its properties that we are going
to use below. The general reference for this section is
\cite{pS94}.

Let $G$ be a locally-finite connected graph. For a path $t$ in $G$
we denote by $Et$ the edge set of $t$. Similarly, by $ET$ we
denote the edge set of a family of paths $T$ in $G$. The
{\emph{extremal length}} of a family of paths $T$ in $G$,
$\lambda(T)$, is defined as
$$
\lambda(T)^{-1}=\inf\bigg\{\sum_{e\in ET}\mu(e)^2\bigg\},
$$
where the infimum is taken with respect to all density functions $\mu$
defined on the edge set $ET$, such that for all $t\in T$,
$$
\sum_{e\in Et}\mu(e)\geq 1.
$$

The extremal length of the family of paths connecting two vertices
or a vertex to infinity, is equal to (a scalar multiple of) the
effective resistance between the two vertices, respectively the
vertex and infinity. It is known that a locally-finite graph $G$
is hyperbolic (parabolic) if and only if $\lambda(T_v)$ is finite
(infinite) for some, and hence every, vertex $v\in VG$, where
$T_v$ is the family of paths connecting $v$ to infinity.

Let $T,\ T_i,\ i\in I$, be families of paths in $G$, where $I$ is
at most countable. We assume that $ET_i\cap ET_j=\emptyset,\ i\neq j$.
Suppose that for every $t\in T$ and every $i\in I$, there exists
$t_i\in T_i$, which is a subpath of $t$. Then (see Appendix~\ref{Proextlen})
\begin{equation}\label{E:Exl}
\lambda(T)\geq\sum_{i\in I}\lambda(T_i).
\end{equation}

\section{Conjecture}

We recall that the conjecture of R. Nevanlinna (\cite{rN70}, p. 312) 
states that
{\emph{a surface $(X, f)$ of the class $F_q$ is of a hyperbolic or
a parabolic type, according to whether the angle geometry of the surface
is ``Lobachevskyan'' or ``Euclidean'', i.e. according to whether the mean
excess $E_m$ is negative or zero}}.

\section{Counterexample 1}\label{S:Nc}\label{S:Ce1}

In what follows, we mean by $a\asymp b$ that there are absolute
positive constants $c_1, c_2$, such that $c_1a\leq b\leq c_2a$;
similarly $a\lesssim b$ means that there is an absolute  positive
constant $c$, such that $ac\leq b$.

\subsection{Speiser Graph}

\begin{figure}[h]
\centerline{\scalebox{1.2}{\includegraphics{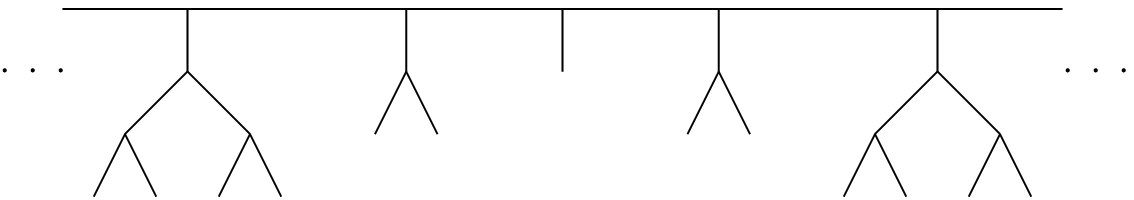}}}
\caption{Graph Tr}
\label{fig:tr}
\end{figure}

First we consider an {\it infinite linear graph}, i.e. an infinite
connected graph where each vertex has degree 2. Next, we fix a
vertex of this graph, and denote it by 0. To a vertex of this
graph that is at a distance $i$ from 0, we attach a binary tree of
$(i+1)$ generations (see Fig.~\ref{fig:tr}). We denote this graph by $Tr$.

\begin{figure}[h]
\centerline{\scalebox{.8}{\includegraphics{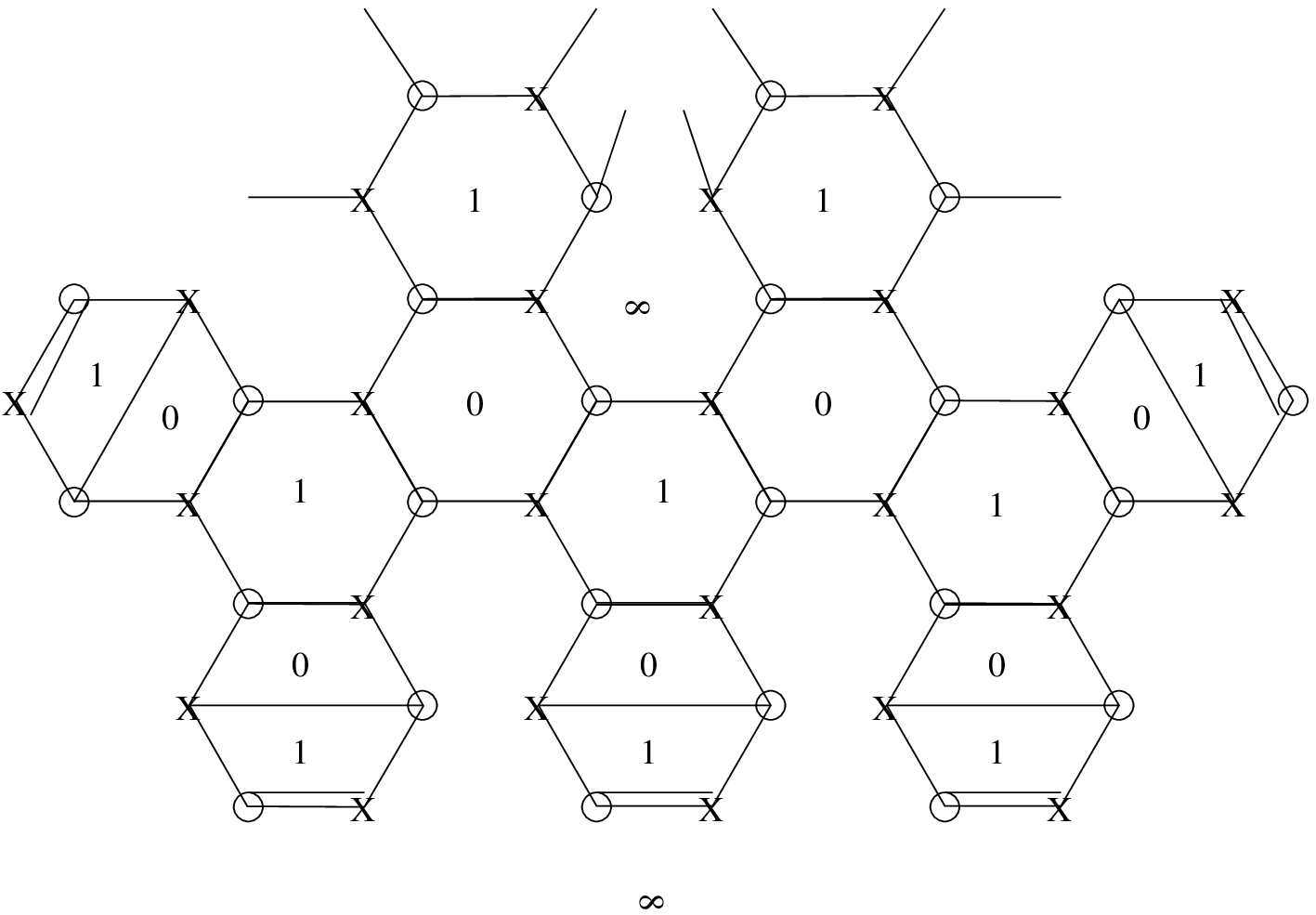}}}
\caption{Graph $\Gamma$}
\label{fig:bee}
\end{figure}

The vertices of $Tr$ of degree one we call {\emph{leaves}}. To
obtain a Speiser graph $\Gamma$ we replace each vertex of $Tr$ by
a hexagon. Adjacent hexagons correspond to the vertices of $Tr$
that are connected by an edge. The hexagons corresponding to
leaves of $Tr$, which we call {\it {free hexagons}}, should be
completed with two edges, to preserve the degree. We add the edges
to each of these hexagons, so that the pair of opposite vertices
of degree 2 is connected by an edge inside the hexagon, and the
remaining vertices of degree 2 are connected by an edge.
The resulting graph $\Gamma$ has degree 3 at all of its vertices
(see Fig.~\ref{fig:bee}).

We label the faces by $0, 1$, and $\infty$. There are exactly two faces
with infinitely many edges, both labelled by $\infty$. Let $(X,
f)$ denote the surface corresponding to $\Gamma$, $(X, f)\in
F_3=F(0, 1, \infty)$; we chose the extended real line as the base
curve. Notice that the function $f$ is analytic. We need to show that the
surface is parabolic and the mean excess is negative.

\subsection{Parabolicity}

\begin{figure}[h]
\centerline{\scalebox{.6}{\includegraphics{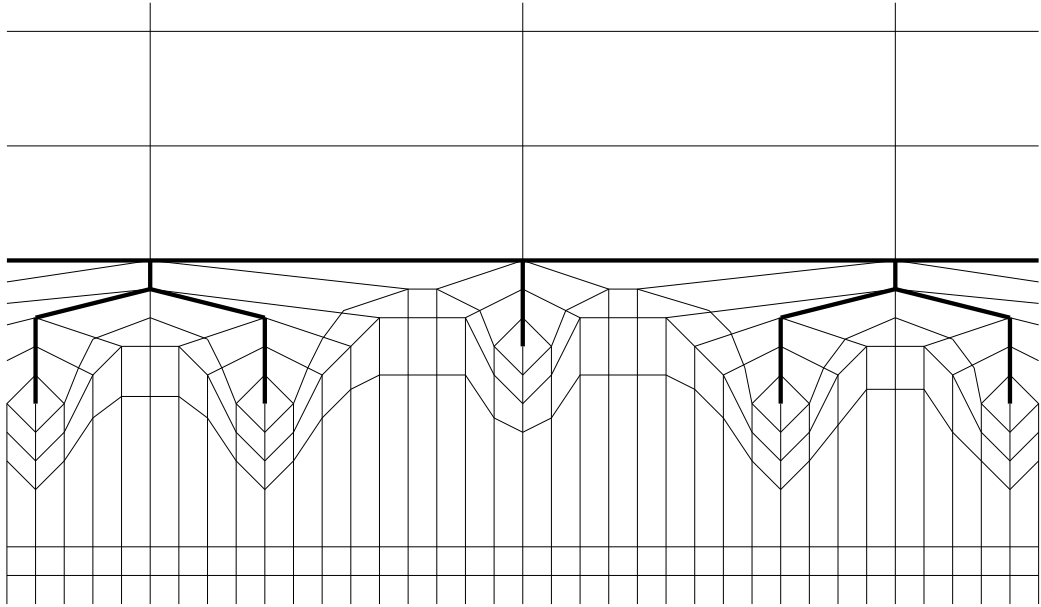}}}
\caption{Dual Graph $\Gamma_4^*$}
\label{fig:ext}
\end{figure}

To prove parabolicity, we make use of Theorem~\ref{T:Ut}. For this
we consider the graph $\Gamma_{4}$. It is easier to deal
with its dual $\Gamma_{4}^*$ (see Fig.~\ref{fig:ext}) though.

\begin{figure}[h]
\centerline{\scalebox{.7}{\includegraphics{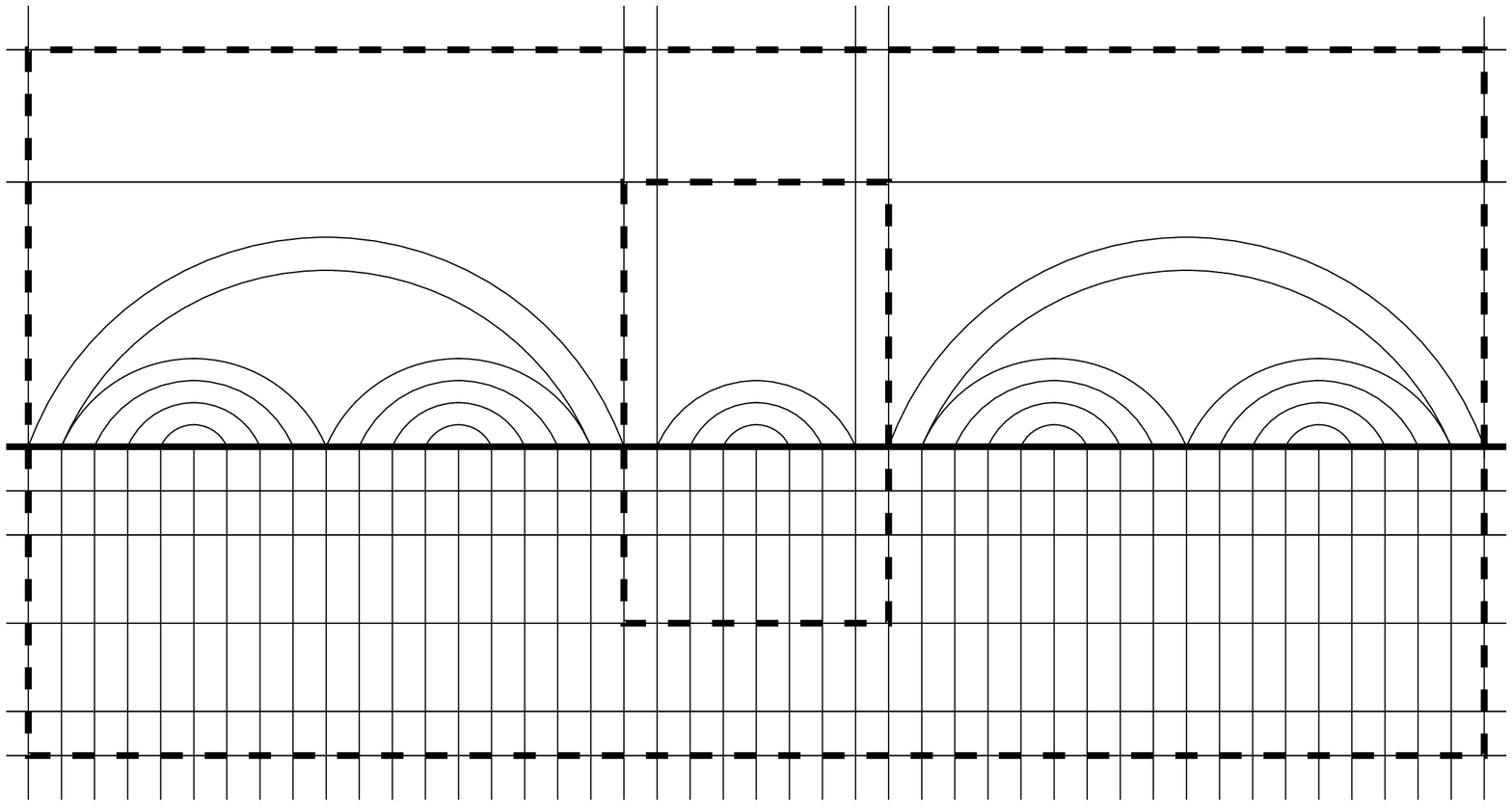}}}
\caption{Graph $\Gamma^*$}
\label{fig:star}
\end{figure}

Since all faces of $\Gamma_{4}$ have a uniformly bounded (by 6)
number of sides, and the degree of $\Gamma_{4}$ is bounded (it is
4), $\Gamma_{4}$ is roughly isometric to $\Gamma_{4}^*$, and hence
they are simultaneously hyperbolic (parabolic). To simplify
further, we pass from $\Gamma_{4}^*$ to a roughly isometric graph
$\Gamma^*$ of bounded degree. The graph $\Gamma^*$ (see Fig.~\ref{fig:star})
consists of a {\emph{coarse}} lattice in the upper half-plane, a
{\emph{fine}} lattice in the lower half-plane, and edges,
which we call {\emph{bridges}}, that connect vertices on the real
line. The bridges are chosen in such a way, that identification of
vertices connected by them gives a rough isometry
$\Gamma^*\to\Gamma_{4}^*$. We denote by $v$ the vertex on the real
line with respect to which $\Gamma^*$ is symmetric, and by $T_v$
the family of paths connecting $v$ to infinity. We show that
$\lambda(T_v)=\infty$, which implies that the surface is
parabolic.

Let $A_i$ be a finite subgraph of $\Gamma^*$, which is an annulus
of combinatorial width 1 in the upper half-plane, of combinatorial
width $\asymp 2^i$ in the lower half-plane, and which contains
bridges. In Figure~\ref{fig:star} the inner and outer boundaries of the first
annulus are marked by dashed lines.

Let $T_i$ denote the family of paths in $A_i$ that connect the
inner and outer boundaries. We consider a density function
$\mu_i$, which assigns the value $1/2^i$ to every edge of $A_i$ in
the lower half-plane, and the value $1$ to every edge in the upper
half-plane. To the bridges we assign values as follows. We say
that a bridge has {\it size} $k$, if it connects the vertices that
are at a distance $k$ with respect to the real line. Now, to a
bridge of size $k$ we assign the value $1/2^{l-1}$, where $l=l(k)$
is the number of bridges of size $k$ in $A_i$. We notice that for
each $l$ there are at most 4 different sizes $k$ for which
$l(k)=l$.

For every path $t_i\in T_i$ we have
$$
\sum_{e\in Et}\mu_i(e)\gtrsim 1.
$$

From the definition of the extremal length we get
\begin{equation}\label{E:El}
\lambda(T_i)^{-1}\leq \sum_{e\in ET_i}\mu_i(e)^2.
\end{equation}
Since there are $\asymp i$ edges of $A_i$ in the upper half-plane and
$\asymp 2^{2i}$ in the lower half-plane, these two parts combined
contribute $\asymp i\times 1+2^{2i}\times (1/2^{2i})\asymp i$ to
the right-hand side of~(\ref{E:El}). The bridges of $A_i$
contribute
$$
\sum_{k}l(k)\frac1{2^{2(l(k)-1)}}\lesssim 1.
$$
Combining the above estimates, we conclude that
$$
\sum_{e\in ET_i}\mu_i(e)^2\asymp i,
$$
and hence
$$
\lambda(T_i)\gtrsim \frac1i.
$$
Therefore $\lambda(T_v)\gtrsim\sum_{i=1}^{\infty}1/i=\infty$.

\subsection{Mean Excess}

Now we show that the (upper) mean excess is negative.

\begin{figure}[h]
\centerline{\scalebox{.7}{\includegraphics{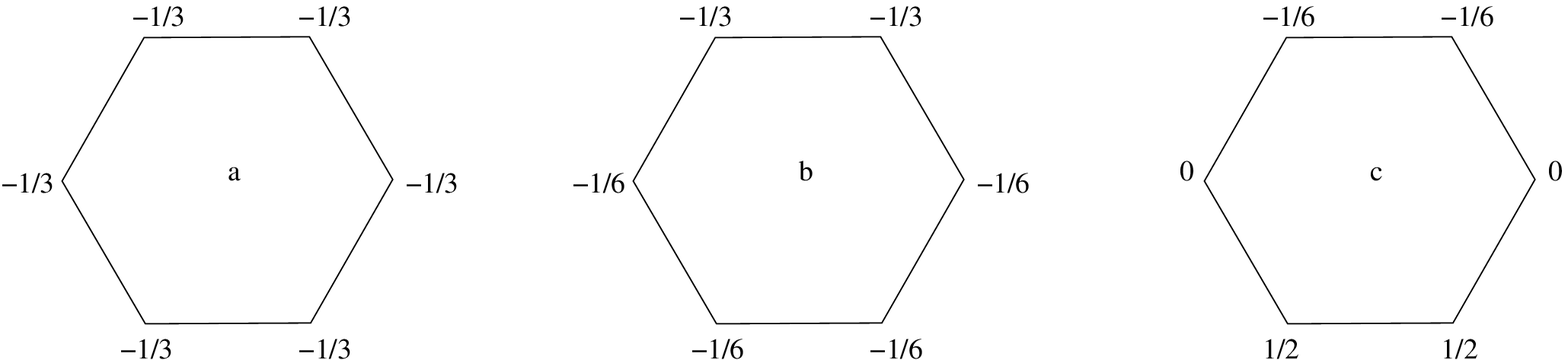}}}
\caption{Hexagons}
\label{fig:hex}
\end{figure}

To each vertex of $Tr$ there corresponds a hexagon of $\Gamma$.
There are 3 types of hexagons, according to the excess assigned to
their vertices. We call these types $a, b$, and $c$ (see Fig.~\ref{fig:hex}, 
where the numbers next to the vertices of hexagons are the corresponding
values of the excess).

\begin{figure}[h]
\centerline{\scalebox{.8}{\includegraphics{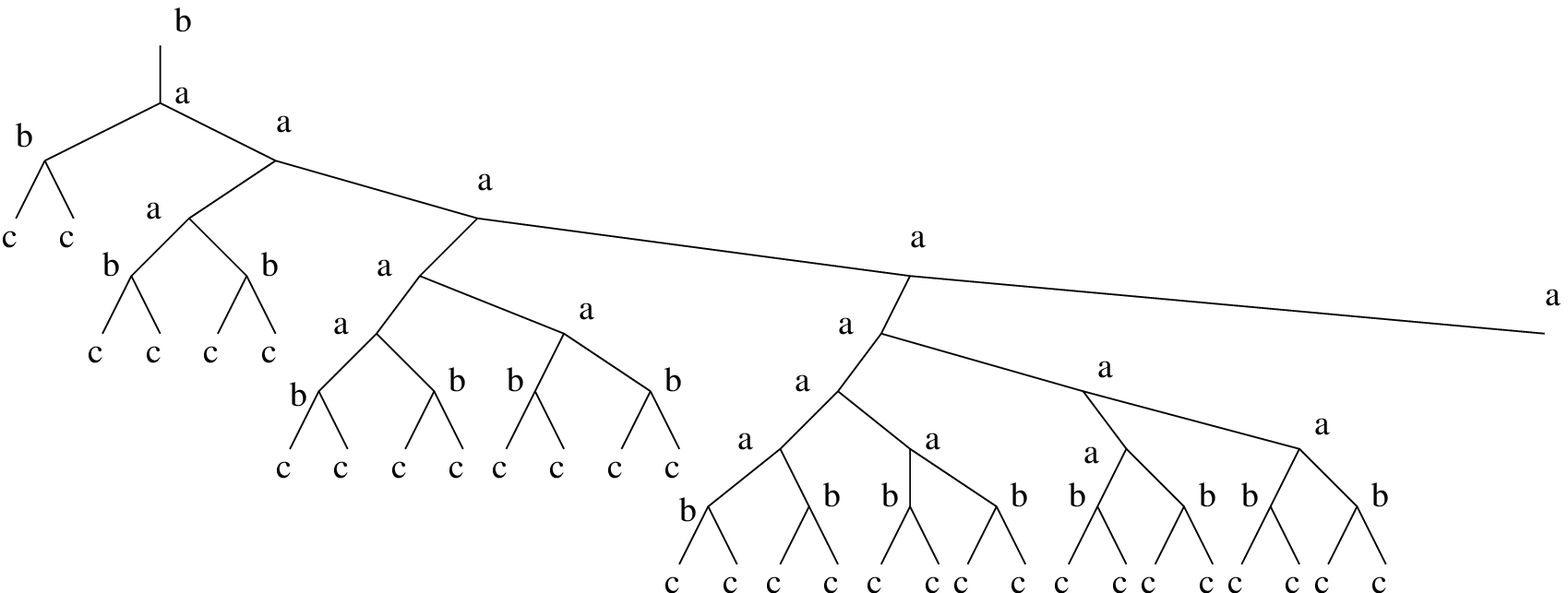}}}
\caption{Tree}
\label{fig:tree}
\end{figure}

In order to compute the mean excess of $\Gamma$, we look at the
graph $Tr$, whose vertices are labelled by $a, b$, and $c$ (due to the
symmetry of $Tr$, we can consider only the part of it which is to
the right of 0, see Fig.~\ref{fig:tree}).

We split the sequence $\Gamma_{(i)}$ of balls into 4 subsequences,
according to the index $\mod 4$, and compute the mean excess in
each case, counting how many hexagons of every type $a, b, c$ are
included:
$$
\begin{aligned}
&1)\ \frac{[4(-\frac13)(2^{k+1}-2)-2(-\frac13)2^{k-1}]+[6(-\frac16)2^{k-1}]
+[2\frac12(2^k-2)]}{[4(2^{k+1}-2)-2(2^{k-1})]+[6(2^{k-1})]
+[4(2^k-2)]}\sim-\frac{11}{84};\\
&2)\ \frac{[4(-\frac13)(2^{k+1}-2)]+[4(-\frac16)2^{k}]
+[2\frac12(2^k-2)]}{[4(2^{k+1}-2)]+[4(2^{k})]
+[4(2^k-2)]}\sim-\frac{7}{48};\\
& 3)\ \frac{[4(-\frac13)(3(2^{k})-2)-2(-\frac13)2^{k}]+[4(-\frac16)2^{k}]
+[2\frac12(2^{k+1}-2-2^k)]}{[4(3(2^{k})-2)-2(2^{k})]+[4(2^{k})]
+[4(2^{k+1}-2)-2(2^k)]}\sim-\frac{3}{20};
\\
&4)\ \frac{[4(-\frac13)(3(2^{k})-2)]+[4(-\frac16)2^{k}]
+[2\frac32(2^{k+1}-2)]}{[4(3(2^{k})-2)]+[4(2^{k})]
+[4(2^{k+1}-2)]}\sim-\frac{1}{9}.
\end{aligned}
$$
Since the limits of all the subsequences are $<0$, the mean excess
$E_m$ is also $<0$.

\section{Counterexample 2: No Asymptotic Values}\label{S:Ce2}

A face of a Speiser graph with infinitely many edges on the boundary 
corresponds to an asymptotic spot of $f$, and a face with finitely
many edges corresponds to a critical point. In the previous example 
we had two asymptotic spots of $f$ with the same value $\infty$. Hence the
function is analytic. We give another example where the function is
meromorphic, but it does not have asymptotic values.

We notice that the graph $\Gamma_4$, which is
an extended Speiser graph of $\Gamma$, is itself a Speiser graph
of degree 4. It provides us with an example of a parabolic surface
$(X, f)\in F_4$, whose mean excess is negative. The new feature is
that $\Gamma_4$ does not have logarithmic elementary regions, and, 
moreover, all algebraic critical points have bounded order.

\section{Counterexample 3}\label{S:construct}

P.~Doyle~\cite{pD84} proved that the surface $(X,\psi)$ is parabolic
if and only if a certain modification of the Speiser graph is recurrent.
(See~\cite{pD84} and~\cite{pS94} for background on recurrence and
transience of infinite graphs.)
In the particular case where $k_f$ is bounded, the recurrence of the
Speiser graph itself is equivalent to $(X,\psi)$ being parabolic.
Though we will not really need this fact,
it is not too hard to see that in a Speiser graph satisfying $E_m<0$
the number of vertices in a ball grows exponentially with the radius.
Thus, we may begin searching for a counterexample by considering
recurrent graphs with exponential growth.
A very simple standard example of this sort is a tree constructed
as follows.
In an infinite $3$-regular tree $T_3$, let $v_0,v_1,\dots$ be
an infinite simple path.  Let $T$ be the set of vertices $u$
in $T_3$ such that $d(u,v_n)\le n$ for all sufficiently large $n$.
Note that there is a unique infinite simple path in $T$ starting from
any vertex $u$.  This implies that $T$ is recurrent.  It is
straightforward to check that the number of vertices
of $T$ in the ball $B(v_0,r)$ grows exponentially with $r$.

\begin{figure}
\tabskip=1em plus2em minus.5em
\halign to\hsize{\hfill#\hfill&\hfill#\hfill\cr
\includegraphics*[height=1.6in]{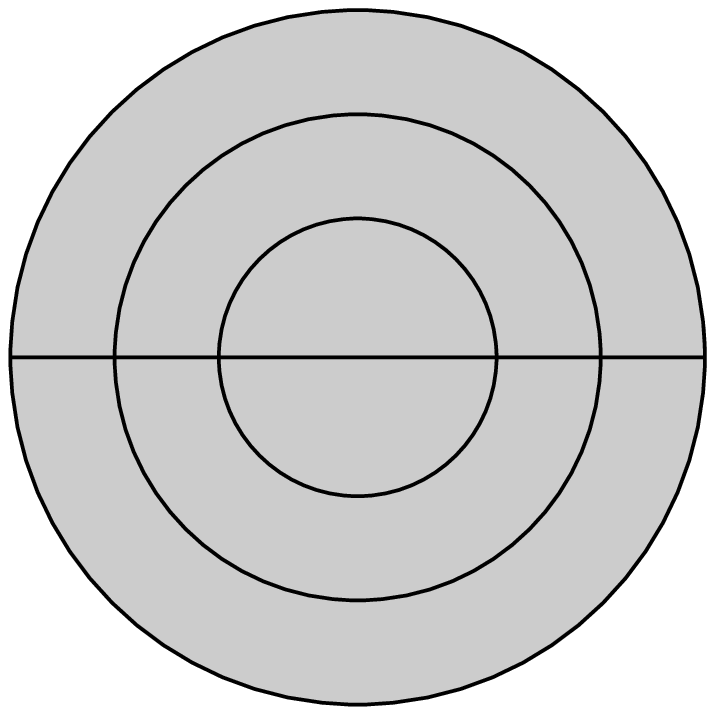}%
&%
\includegraphics*[height=1.6in]{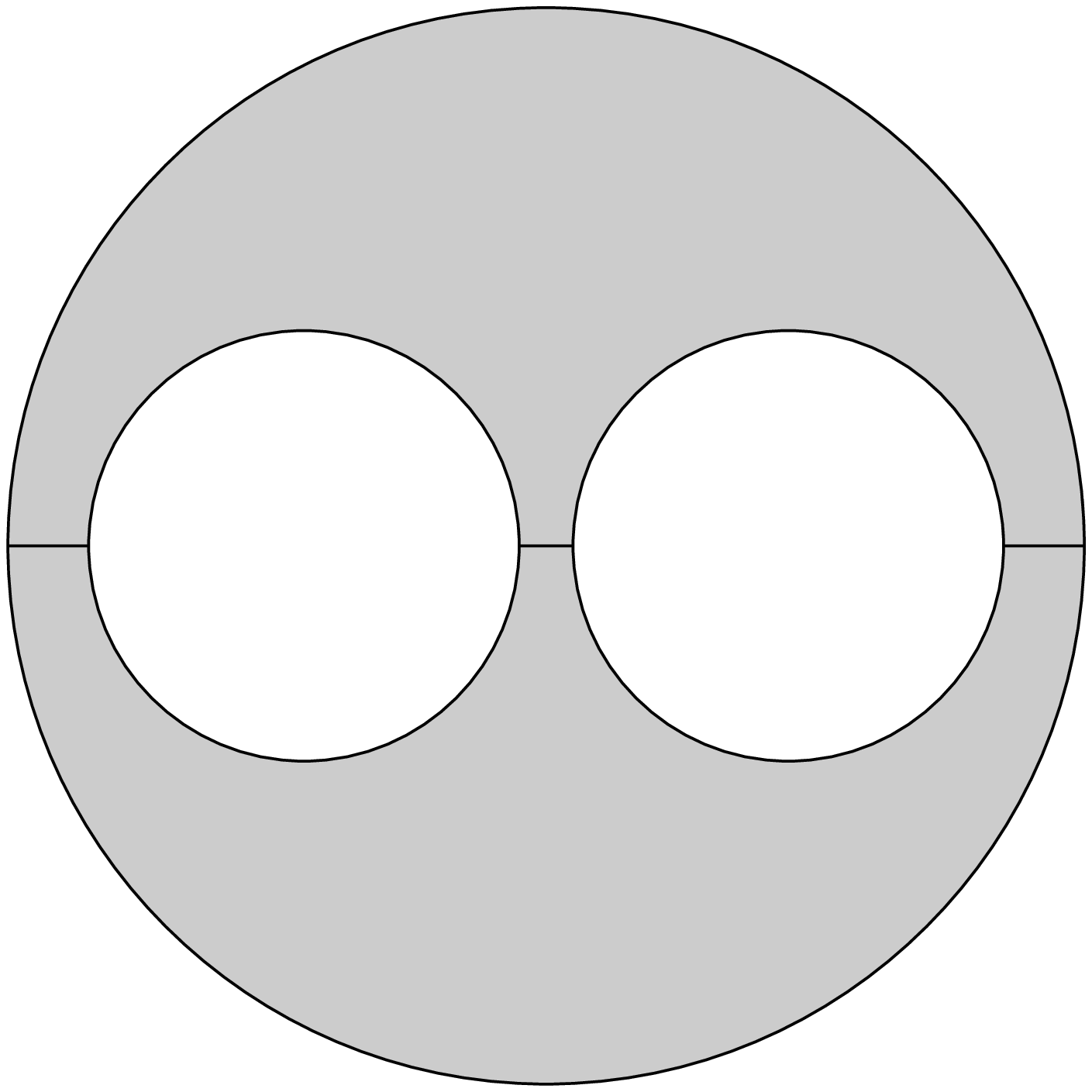}%
\cr
(a)&(b)\cr
}
\caption{\label{f.pantsandfinger}The surfaces $S(v)$.
}
\end{figure}

\begin{figure}
\SetLabels
(.29*.58)$S(v_1)$\\
(.4*.7)$S(v_2)$\\
(.65*.9)$S(v_3)$\\
\endSetLabels
\centerline{\AffixLabels{%
\includegraphics*[height=3in]{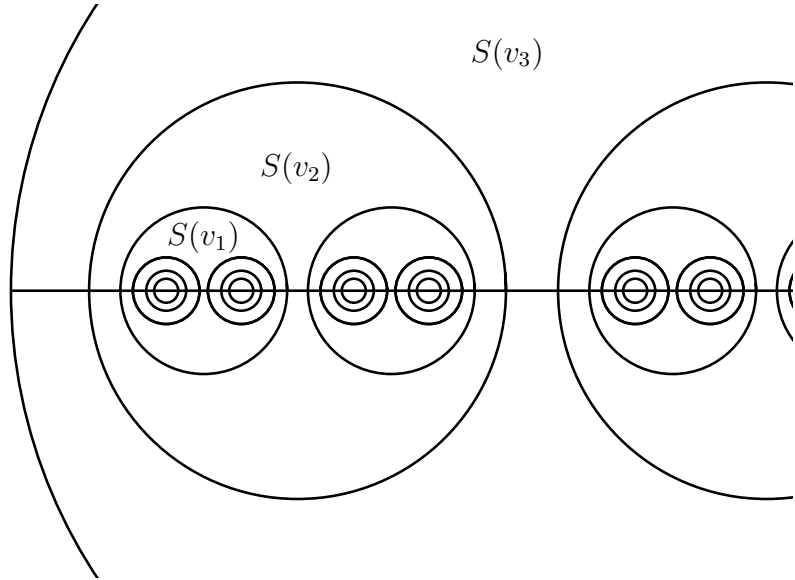}
}}
\caption{\label{f.Speiser}The Speiser graph with $s=2$.
}
\end{figure}

Our Speiser graph counterexample
is a simple construction based on the tree $T$.
Fix a parameter $s\in\{1,2,\dots\}$, whose choice will be discussed
later.  To every leaf (degree one vertex) $v$ of $T$ associate a closed
disk $S(v)$ and on it draw the graph indicated in
Figure~\ref{f.pantsandfinger}.(a),
where the number of concentric circles, excluding
$\partial S(v)$, is $s$.
If $v$ is not a leaf, then it has degree $3$.
We then associate to it the graph indicated in
Figure~\ref{f.pantsandfinger}.(b),
drawn on a triply connected domain $S(v)$.  We combine these to form
the Speiser graph $\Gamma$ as indicated in figure~\ref{f.Speiser},
by pasting the outer boundary of the surface corresponding to each vertex
into the appropriate inner boundary component of its parent.
Here, the parent of $v$ is the vertex $v'$ such that
$d(v',v_n)=d(v,v_n)-1$ for all sufficiently large $n$.

Every vertex of $\Gamma$ has degree $4$ and every face
has $2$, $4$, or $6$ edges on its boundary.  Therefore,
$\Gamma$ is a Speiser graph.
Consequently, as discussed above, there is a surface
spread over the sphere $X=(\R^2,\psi)$ whose Speiser graph
is $\Gamma$.
It is immediate to verify that $\Gamma$ is recurrent,
for example, by the Nash-Williams criterion.
Doyle's Theorem~\cite{pD84} then implies that $X$ is parabolic.
Alternatively, one can arrive at the same conclusion
by noting that there is an infinite sequence of
disjoint isomorphic annuli on $(\R^2,\Gamma)$ separating
any fixed point from $\infty$, and applying extremal length.
(See~\cite{lA73},~\cite{oL73} for the basic properties of extremal
length.)

We now show that $\overline E_m<0$ for $\Gamma$.  Note that
the excess is positive only on vertices on the boundaries
of $2$-gons, which arise from leaves in $T$.
On the other hand, every vertex of degree $3$ in $T$ gives
rise to vertices in $\Gamma$ with negative excess.  Take
as a basepoint for $\Gamma$ a vertex $w_0\in S(v_1)$ with negative excess. 
It is easy to see that there are constants
$a>1,c>0$, such that the number $n^-_r$ of negative excess vertices
in the combinatorial ball $B(w_0,r)$ about $w_0$ satisfies
$c\,a^r\le n^-_r\le a^r/c$.

If $w$ is a vertex with positive excess, then there is a unique
vertex $\sigma(w)$ with negative excess closest to $v$;
in fact, if $w\in S(v)$, then $\sigma(w)$ is the closest
vertex to $w$ on $\partial S(v)$,
and the (combinatorial) distance from $w$ to $\sigma(w)$ is our parameter
$s$.  The map $w\mapsto \sigma(w)$ is clearly injective.
This implies that the number $n^+_r$ of positive excess
vertices in $B(w_0,r)$ satisfies
$n^+_{r+s}\le n^-_r$, $r\in\{0,1,2,\dots\}$.
By choosing $s$ sufficiently large, we may therefore
arrange to have the total excess in $B(w_0,r)$
to be less than $-\epsilon\, a^r$, for some $\epsilon>0$
and every $r\in\{0,1,2,\dots\}$.
It is clear that the number of vertices with zero excess in $B(w_0,r)$
is bounded by a constant (which may depend on $s$)
times $n^-_r$.  Hence $\overline E_m<0$ for $\Gamma$.

By allowing $s$ to depend on the vertex in $T$, if necessary,
we may arrange to have $\underline E_m=\overline E_m$;
that is, $E_m$ exists, while maintaining $E_m<0$.
We have thus demonstrated that the resulting surface
is a counterexample in $F_4$ to the second implication in
Nevanlinna's problem.

\section{A Non-Positive Curvature Example}\label{S:nonpos}

We now construct an example of a simply connected,
complete, parabolic surface $Y$ of nowhere positive
curvature, with the property
\begin{equation}\label{E:Intc}
\int_{D(a,r)}\text{curvature} <-\epsilon\,\text{area}\bigl(D(a, r)\bigr),
\end{equation}
for some fixed $a\in Y$ and every $r>0$, where $D(a, r)$ denotes the
open disc centered at $a$ of radius $r$, and $\epsilon>0$
is some fixed constant.

Consider the surface $\C=\R^2$ with the metric $|dz|/y$ in
$P=\{z=x+iy:\ y\geq1\}$,
and $\exp{(1-y)}|dz|$ in $Q=\{y<1\}$. We denote this surface by $Y$.
Let $\beta$ denote
the curve $\{y=1\}$ in $Y$, i.e., the common boundary of $P$ and $Q$.

Let $Q'$ denote the universal cover of $\{z\in\C:|z|>1\}$.
Note that $Q$ is isometric to $Q'$ via the map $z\mapsto \exp(iz+1)$.
Hence the curvature is zero on $Q$, and the geodesic curvature
of $\partial Q$ is $-1$.  The geodesic curvature of $\partial P$ is $1$.
Consequently, $Y$ has no concentrated curvature on $\beta$.
The surface $Y$ is thus a ``surface of bounded curvature'',
also known as an Aleksandrov
surface (see~\cite{aA67}, \cite{yR93}).
The  curvature measure of $Y$ is absolutely
continuous with respect to area; the curvature of $Y$ is -1 (times the area
measure) on $P$ and $0$ on $Q$.

The surface $Y$ is parabolic,
and the uniformizing map is the identity map onto $\R^2$ with the standard
metric.

We will now prove (\ref{E:Intc}) with $a=i$.
Set $\beta_r=D(a,r)\cap\beta$. Note that the shortest path in $Y$ between
any two points on $\beta$ is contained in $P$, and is the arc
of a circle orthogonal to $\{y=0\}$. Using the Poincare disc model,
it is easy to see that there exists a constant $c>0$, such that
\begin{equation}\label{lbr}
c\,e^{r/2}\leq \text{length}\beta_r\leq e^{r/2}/c,
\end{equation}
where the right inequality holds for all $r$, and the left for all sufficiently
large $r$.
By considering the intersection of $D(a,r)$ with the strip
$1<y<2$ it is clear that
\begin{equation}\label{E:Plb}
O(1)\text{area}\bigl(P\cap D(a,r)\bigr)\ge \text{length}\beta_r,
\end{equation}
for all sufficiently large $r$.

Consider some point $p\in Q$, and let
$p'$ be the point on $\beta$ closest to $p$.
It follows easily (for example, by using the
isometry of $Q$ and $Q'$) that if $q$ is any point in $\beta$,
then $d_Q(p,q)=d_Q(p,p')+d_Q(p',q)+O(1)$.
Consequently, if $d(p,a)\le r$,
then there is an $s\in[0,r]$ such that
$p'\in \beta_s$ and $d_Q(p,p')\le r-s+O(1)$.
Furthermore,
it is clear that the set of points $p$ in $Q$ such
that $p'\in \beta_s$ and $d_Q(p,p')\le t$ has
area $O(t^2+t)\,\text{length}\beta_s$.  Consequently,
$$
\text{area}\bigl(Q\cap D(a,r)\bigr)
\le O(1) \sum_{j=0}^{r} (j+1)^2\text{length}\beta_{r-j}\,.
$$
Using~(\ref{lbr}), we have
\begin{equation}\label{E:Qub}
\text{area}\bigl(Q\cap D(a,r)\bigr)\le O(1)\,\text{length}\beta_r,
\end{equation}
for all sufficiently large $r$.

Now, combining~(\ref{E:Plb}) and (\ref{E:Qub}),
we obtain (\ref{E:Intc}) for all sufficiently large $r$. It therefore
holds for all $r$.

\newpage
\vspace*{2cm}
\begin{center}
This page deliberately left blank
\end{center}

\chapter{Analytic Endomorphisms}

In this chapter we study the question of recovering a domain from its
semigroup of analytic endomorphisms.


\section{Semigroups}

If $\Omega$ is a domain in $\C^n$, its \emph{analytic endomorphism}
is an analytic map from $\Omega$ into itself. Analytic endomorphisms
of $\Omega$ form a semigroup (with the identity map as the unit),
which we denote by
$E(\Omega)$. It is non-comutative. An isomorphism between semigroups is
a map that preserves the operation and sends the unit to the unit.

A \emph{biholomorphic} map between two domains is a one-to-one onto analytic
map, whose inverse is also analytic. We say that a map is
\emph{antibiholomorphic}, if its complex conjugate is biholomorphic.
By an \emph{(anti-)biholomorphic} map we mean a map which is either
biholomorphic or antibiholomorphic. It is obvious that if
$\psi:\ \Omega_1\to\Omega_2$ is an (anti-)biholomorphic map,
then the map between the corresponding semigroups
$f\mapsto\psi\circ f\circ\psi^{-1}$ is an isomorphism.

We study the converse implication, namely given an isomorphism between
semigroups, is there an (anti-)biholomorphic map that conjugates it?

\section{Statement of the Theorem}

\begin{theorem}\label{T:Mt}
Let $\Omega_1,\ \Omega_2$ be bounded domains in $\C^n,\ \C^m$
respectively, and suppose that there exists
$\varphi:\ E(\Omega_1)\rightarrow E(\Omega_2)$,  an isomorphism of
semigroups. Then $n=m$ and there exists an (anti-)biholomorphic
map $\psi:\ \Omega_1\rightarrow\Omega_2$ such that
\begin{equation}\label{E:Con}
\varphi f=\psi\circ f\circ \psi^{-1}, \ \ \text{for all}\ f\in E(\Omega_1).
\end{equation}
\end{theorem}

\section{Topology}\label{S:Top}

\subsection{Constant Endomorphisms}

To construct the (anti-)biholomorphic map and deduce the desired properties
of it, we need to express certain properties of elements of the semigroup,
such as an element being
injective, or constant, in terms of the semigroup structure. This will allow
us to conclude that an element with a property, say being injective,
will map to an element with the same property. The most crucial property
for a construction of the biholomorphic map is constantness of an element.

For a bounded domain $\Omega$ in $\C^n$, we denote by $C(\Omega)$ the
subsemigroup of $E(\Omega)$ consisting of constant maps. An endomorphism
$c_z$ is constant if it sends $\Omega$ to a point $z\in\Omega$. The subset
$C(\Omega)\subset E(\Omega)$ can be described using only the semigroup
structure as follows:

\begin{equation}\label{E:Co}
c\in C(\Omega) \ \text{iff} \ \ \forall (f\in E(\Omega)), \ \ (c\circ f=c).
\end{equation}
In other words, a constant map is a left zero.

It is clear that we have a bijection between constant endomorphisms of
$\Omega$ and points of this domain as a set: to each $z$ corresponds a
unique $c_z\in C(\Omega)$ and vice versa, so we can identify the two.
Under this identification, a subset of $\Omega$ corresponds to a
subsemigroup of $C(\Omega)$.

\subsection{Construction of $\psi$}

Having defined points of a domain in terms of its semigroup structure of
analytic endomorphisms, we can construct a map $\psi$ between $\Omega_1$
and $\Omega_2$ as follows:

\begin{equation}\label{E:Dp}
\psi(z)= w \ \ \text{iff}\ \ \varphi c_z=c_w.
\end{equation}

So defined, $\psi$ satisfies (\ref{E:Con}). Indeed, let $f\in E(\Omega_1),\,
f(z)=\zeta
$. This is equivalent to
\begin{equation}\label{E:Se}
f\circ c_z=c_{\zeta}.
\end{equation}
Applying $\varphi$ to both sides of (\ref{E:Se}) we have

\begin{equation}\label{E:3}
\varphi f\circ c_{\psi(z)}=c_{\psi(\zeta)}.
\end{equation}
But (\ref{E:3}) is equivalent to $\varphi f(\psi(z))=\psi(\zeta)=\psi(f(z))$,
which is (\ref{E:Con}).

\subsection{Continuity of $\psi$}

We describe the topology of a domain $\Omega$ using its injective
endomorphisms. A map $f\in E(\Omega)$ is injective if and only if
$$
\forall (c'\in C(\Omega)),\ \forall (c''\in C(\Omega)),\ \
((f\circ c'= f\circ c'')\Rightarrow(c'= c'')).
$$
We denote the class of injective endomorphisms of $\Omega$ by
$E_i(\Omega)$. For every $f\in E_i(\Omega),\ f_i(\Omega)$ is open
\cite{sB48}.
The family
$
\{f(\Omega), \ \ f\in E_i(\Omega)\}
$
of subsets of $\Omega$
forms a base of topology, because every $z\in \Omega$ has a neighborhood
$f(\Omega)$, where $f(\zeta)=z + \lambda(\zeta-z)$, $f$ belongs to
$E_i(\Omega)$ for every $\lambda$ such that $|\lambda|$ is small. This is
the place where we use the boundedness of domains.

Thus we described subsets of $\Omega$ and the topology on it
using only the semigroup structure of $E(\Omega)$. Since this is so, the
semigroup structure also defines the notions of an open set, closed set,
compact set, and closure of a set.

Now we can easily prove continuity of the map $\psi$ constructed above.
Indeed, let $g(\Omega_2), \ g\in E_i(\Omega_2)$ be a set from the base of
topology of $\Omega_2$. We take $f=\varphi^{-1}g$. Then
$f\in E_i(\Omega_1)$ and $\psi^{-1}(g(\Omega_2))=f(\Omega_1)$, which
proves that $\psi$ is continuous. Since $\varphi$ is an isomorphism,
the same argument works to prove that $\psi^{-1}$ is also continuous,
and thus $\psi$ is a homeomorphism.

Therefore the domains $\Omega_1, \ \Omega_2$ are homeomorphic, and hence
\cite{wH41} they have the same dimension, i.e. $n = m$.

\section{Localization}\label{S:Loc}

In order to prove that $\psi$ is (anti-)biholomorphic, we will introduce
a system of projections. The main difficulty in extracting any useful
information from such a system is that a projection in general does not
have to be an endomorphism. To overcome this difficulty, the
following localization lemma proves useful.

\begin{lemma}\label{L:E}
Suppose $H$ is a semigroup with identity, and $f$ an element of $H$ with
the following two properties:

(i) $hf=fh$, for every $h$ in $H$, and

(ii) $h_1f=h_2f$ implies $h_1=h_2$, for every $h_1$ and $h_2$ in $H$.

Then there exists a semigroup $S_f$ and a monomorphism $i:\ H\rightarrow
S_f$, such that $i(f)$ is invertible in $S_f$ and commutes with all
elements of
$S_f$. Moreover, the semigroup $S_f$ satisfies the following universal
property: for every semigroup $S_1$ with a monomorphism
$i_1: \ H\rightarrow S_1$ such that $i_1(f)$ is invertible in $S_1$
and commutes with all elements of $S_1$, there exists a unique
monomorphism $\hat i_1: \ S_f\rightarrow S_1$ such that
$i_1=\hat i_1\circ i$.
\end{lemma}

\begin{remark}
Uniqueness of $\,\hat i_1$ implies that the semigroup $S_f$ with the
universal property is unique up to isomorphism.
\end{remark}

\emph{Proof.}
We construct $S_f$ as follows. First we consider formal expressions of the
form $hf^k$, where $h\in H$ and $k$ is an integer (may be positive, negative
or
zero). Then we define a multiplication on this set:
$h_1f^{k_1}*h_2f^{k_2}=h_1h_2f^{k_1+k_2}$. Next we consider a relation
on the set of formal expressions: $h_1f^{k_1}\sim h_2f^{k_2}$ if
$k_1\leq k_2$ and $h_1=h_2f^{k_2-k_1}$ in $H$, or $k_2\leq k_1$ and
$h_2=h_1f^{k_1-k_2}$ in $H$. It is easy to verify that this is an
equivalence relation and it is compatible with the operation $*$;
that is, $x\sim y, \ u\sim v$ implies $x*u\sim y*v$.

Lastly, let $S_f$ be the set of equivalence classes with the binary
operation induced by $*$. For $S_f$ to be a semigroup, we need to show
that the binary operation $*$ is associative. Let
$h_1f^{k_1}\sim h_1'f^{k_1'}$, $h_2f^{k_2}\sim h_2'f^{k_2'}$, and
$h_3f^{k_3}\sim h_3'f^{k_3'}$. We need to show that
$(h_1f^{k_1}*h_2f^{k_2})*h_3f^{k_3}\sim h_1'f^{k_1'}*
(h_2'f^{k_2'}*h_3'f^{k_3'})$. By the definition of the operation $*$,
the last equivalence is the same as $h_1h_2h_3f^{k_1+k_2+k_3}
\sim h_1'h_2'h_3'f^{k_1'+k_2'+k_3'}$. Assuming that $k_1+k_2+k_3
\leq k_1'+k_2'+k_3'$, we have essentially one possibility to consider
(the others are either similar or trivial): $k_1\leq k_1', \ k_2\leq
k_2'$, and $k_3'\leq k_3$. In this case, $h_1h_2h_3f^{k_3-k_3'}=
h_1'h_2'h_3'f^{k_1'-k_1+k_2'-k_2}$. Now we can use the cancellation
property (ii) to get the desired equivalence.

The semigroup $H$ is embedded into $S_f$ via $i: \ h\mapsto [hf^0]$.
The element $i(f)=[\text{id} f]$, where $\text{id}$ is the identity in
$H$, is invertible in $S_f$ with the inverse $[\text{id} f^{-1}]$.
Clearly, $[{\text{id}}f]$ commutes with all elements of $S_f$.

Now suppose that $S_1$, $i_1:\ H\rightarrow S_1$, is a semigroup and a
monomorphism, such that $i_1(f)$ is invertible in $S_1$ and commutes
with all elements of $S_1$. Then we define
$$
\hat i_1([hf^k])=i_1(h)(i_1(f))^k.
$$
This definition does not depend on a representative of $[hf^k]$. Indeed,
suppose $h_1f^{k_1}\sim h_2f^{k_2}$ and assume $k_1\leq k_2$.
Then $h_1=h_2f^{k_2-k_1}$, and thus $i_1(h_1)=i_1(h_2)i_1(f)^{k_2-k_1}$.
Hence $i_1(h_1)i_1(f)^{k_1}=i_1(h_2)i_1(f)^{k_2}$.

So defined, $\hat i_1$ is a homomorphism:

\begin{align}
\hat i_1
&([h_1f^{k_1}][h_2f^{k_2}])=\hat i_1([h_1h_2f^{k_1+k_2}])\notag \\
&=i_1(h_1h_2)i_1(f)^{k_1+k_2}=i_1(h_1)i_1(h_2)i_1(f)^{k_1}i_1(f)^{k_2}\notag\\
&=i_1(h_1)i_1(f)^{k_1}i_1(h_2)i_1(f)^{k_2}=\hat i_1([h_1f^{k_1}])\hat i_1
([h_2f^{k_2}]).\notag
\end{align}

The relation $\hat i_1\circ i= i_1$ holds, since $\hat i_1([hf^0])=i_1(h)$
for all $h\in H$.

Uniqueness of $\hat i_1$ is clear, and
Lemma \ref{L:E} is proved. $Box$

We are going to apply this lemma to get an extension of the isomorphism
$\varphi$ restricted to the commutant of an element $f$ to a larger
semigroup that would contain a system of projections.

\section{Extending $\varphi$}\label{S:Ext}

\subsection{Good Elements}

Here we introduce a subsemigroup, whose elements, following \cite{aE93},
we call to be `good'. They are termed `good' because, first of all,
their analytic properties will be useful for us when extending the restricted
isomorphism $\varphi$, and second, all these properties can be expressed
in terms of the semigroup structure.

We say that for a bounded domain $\Omega$ an
element $f\in E(\Omega)$ is {\it{good at}} $z\in \Omega$, denoted by
$f\in G_z(\Omega)$, if

\begin{enumerate}
\item $z$ is a unique fixed point of $f$,
\item $f(\Omega)$ has compact closure in $\Omega$, and
\item $f$ is injective in $\Omega$.
\end{enumerate}

Property 3 of a good element was already stated in terms of the semigroup
structure. Since the topology on $\Omega$ was described using
only the semigroup structure, Property 2 can also be stated in these
terms. Property 1 can be expressed in terms of the semigroup structure as
$$
(f\circ c_z=c_z)\wedge((f\circ c_{\zeta}=c_{\zeta})\Rightarrow
(c_{\zeta}=c_z)).
$$

Since $f$ is an endomorphism of a domain, all eigenvalues $\lambda$ of its
linear part at $z$ satisfy $|\lambda|\leq 1$ \cite{sK98}. Moreover,
$|\lambda|< 1$
because the closure of $f(\Omega)$ is a compact set in $\Omega$. The
injectivity of $f$ implies \cite{sB48} that it is biholomorphic onto
$f(\Omega)$ and the Jacobian determinant of $f$ does not vanish at any
point of $\Omega$.

It is clear that for every $z\in \Omega$ a good element $f$ at $z$ exists.
For example, we can take $f(\zeta)=z+\lambda (\zeta-z)$ with sufficiently
small $|\lambda|$.

\subsection{Extending a Comutant}

Consider a good element $f\in G_z(\Omega)$ and its commutant $H_f(\Omega)$
in $E(\Omega)$:
$$
H_f(\Omega)=\{h\in E(\Omega):\ \ hf=fh\}.
$$
Clearly $H_f(\Omega)$ is a subsemigroup of $E(\Omega)$. The element $f$,
being good (hence injective), satisfies the cancellation property $(ii)$
of Lemma~\ref{L:E} in $H_f(\Omega)$. Thus, by Lemma~\ref{L:E},
we have the extension
$S_f$
of $H_f(\Omega)$ in which $f$ is invertible and commutes with all elements
of $S_f$. In the case of analytic endomorphisms we can embed $H_f(\Omega)$
into the subsemigroup of $A_z$, the semigroup of germs of analytic mappings
at $z$ under composition, consisting of elements that commute with the germ
of $f$ and containing the
germ of $f^{-1}$. We use the universal property of Lemma~\ref{L:E}
to conclude
that $S_f$ is isomorphic to a subsemigroup of $A_z$. We identify $S_f$ with
this semigroup, i.e. we consider elements of $S_f$ as germs of
analytic mappings at $z$.

\subsection{Extending the Isomorphism}

In proving that $\psi$ is (anti-) biholomorphic we need to show that it is
so in a neighborhood of every point of $\Omega_1$. Since an (anti-)
biholomorphic type of a domain is preserved by translations in $\C^n$,
it is enough to show that $\psi$ is (anti-) biholomorphic in a neighborhood
of $0\in \C^n$, assuming that $\Omega_1$ and $\Omega_2$ contain 0 and
$\psi(0)=0$.

Let $\varphi:\ E(\Omega_1)\rightarrow E(\Omega_2)$ be an isomorphism of the
semigroups, $f$ a good element, $f\in G_0(\Omega_1)$, and $H_f(\Omega_1)$
the commutant of $f$. Then clearly $H_g(\Omega_2)=\varphi(H_f(\Omega_1))$
is the commutant of $g=\varphi f$.
By Lemma~\ref{L:E}, we have the extensions $S_f, \ S_g$ of $H_f(\Omega_1)$ and
$H_g(\Omega_2)$ respectively, and by the universal property of this lemma
the isomorphism $\varphi$ extends to an isomorphism
$$
\Phi:\ \ S_f\rightarrow S_g.
$$

\section{System of Projections and Linearization}\label{S:Sys}

\subsection{Very Good Elements}

Let $\Omega$ be a bounded domain in $\C^n$. We say that a good  element
$f\in G_0(\Omega)$ is {\it{very good at}} 0, and write $f\in VG_0(\Omega)$,
if the corresponding semigroup $S_f\subset A_0$ constructed in
Section~\ref{S:Ext}
contains a system of elements, which we call a system of projections,
$\{p_i\}_{i=1}^{n}$ with the following properties:

(a) $\forall \ (i=1,\dots, n),\ \ p_i\neq 0$,

(b) $\forall \ (i=1,\dots, n),  \ \ p_i^2=p_i$, and

(c) $\forall\ (i,\ j=1,\dots, n,\ i\neq j), \ \ p_ip_j=0$.

There does exist a very good element, since we can take $f$ to be a
homothetic transformation at 0 with sufficiently small coefficient, and
$p_i$ a
projection on the $i$'th coordinate of the standard coordinate system.
Clearly, $p_i f = f p_i$ and there exists $k$ such that
$p_if^k\in E(\Omega)$, and hence $p_i\in S_f$. From now on, we fix a very
good element $f\in VG_0(\Omega)$, associated semigroups $H_f(\Omega),\ S_f$
and a system of projections $\{p_i\}$.

We introduce another subsemigroup of $E(\Omega)$:
$$
P_f(\Omega)=\{h\in G_0(\Omega)\cap H_f(\Omega),\ \ hp_i=p_ih,\ \  i=1,\dots,
n\},
$$
where the commutativity relations are in $S_f\subset A_0$. Notice that
$P_f(\Omega)\neq\emptyset$ since $f$ belongs to it.

\subsection{Linearization Lemma}

\begin{lemma}\label{L:Lin}
For every $h\in P_f(\Omega)$ there exists a biholomorphic germ $\theta_h$ at
$0\in \C^n$ such that $\theta_h h=\Lambda\theta_h$, where $\Lambda=
{\rm{diag}}(\lambda_1,\dots, \lambda_n)$ is an invertible diagonal matrix
which is similar to $dh(0)$ in $GL(n,\C)$.
\end{lemma}

\emph{Proof.}
The relations $p_i\neq0,\ p_i^2=p_i$, and $p_ip_j=0, \ i\neq j$,
imply that for
$P_i=dp_i(0)$, the linear part of $p_i$ at 0, we have
$P_i\neq 0, \ P_i^2=P_i$, and $P_iP_j=0, \ i\neq j$. Since the matrices $P_i$
commute, there exists \cite{kH71} a matrix $A\in GL(n, \C)$ such that
$P_i'=AP_iA^{-1}=\Delta_i=\text{diag}(0,\dots, 1,\dots, 0)$, where the only
non-zero entry appears in the $i$'th place.

Since $p_i^2=p_i,\ i=1,\dots, n$, we can use the argument given in
\cite{sK98} to linearize $p_i$, i.e. there exists a biholomorphic
germ $\xi_i$ at 0 such that $\xi_ip_i=P_i\xi_i,\ d\xi_i(0)=\text{id}, \
i=1,\dots, n$. The map $\xi_i$ is constructed in \cite{sK98} as follows:
$$
\xi_i=\text{id}+(2P_i-\text{id})(p_i-P_i), \ \ i=1,\dots, n.
$$
If we take $\xi_i'=A\xi_i$, we have $\xi_i'p_i=P_i'\xi_i'$.
For simplicity of notations, we assume that $\xi_i$ itself conjugates $p_i$
to a diagonal matrix, that is, $P_i=P_i'$ (in this case $P_i$ is not
necessarily $dp_i(0)$, but rather $Adp_i(0)A^{-1},\ d\xi_i(0)=A$).
For every $i=1,\dots, n$, we have $h_iP_i=P_ih_i$, where
$h_i=\xi_ih\xi_i^{-1}$. Let $H_i=dh_i(0)$. Then $H_iP_i=P_iH_i$, and hence
in the $i$'th row and the $i$'th column the matrix $H_i$ has only one
non-zero entry, $\lambda_i$, which is located at their intersection.
Thus $\lambda_i$ has to be an eigenvalue of $H_i$, and hence of the linear
part of $h$. In particular, $0<|\lambda_i|<1$.

Let $I_i:\ \C\rightarrow\C^n$ be the embedding $z\mapsto
(0,\dots, z, \dots, 0)$, where the only non-zero entry is $z$, which is
in the $i$'th place; and $\pi_i:\ \C^n\rightarrow \C$, a projection
$(z_1,\dots, z_n) \mapsto z_i$, corresponding to the $i$'th axis. For
every $i=1,\dots, n$, the map $\pi_ih_iI_i$ sends a neighborhood of 0 in
$\C$ into $\C$, and its derivative at 0, $ \lambda_i$, is an
eigenvalue of $h$. Hence (\cite{lC93}, p. 31) $\pi_ih_iI_i$ is linearized
by the
unique solution $\eta_{h,i}$ of the Schr\"oder equation
\begin{equation}\label{E:4}
\eta(\pi_ih_iI_i)=\lambda_i\eta, \ \ \eta(0)=0,\ \eta'(0)=1.
\end{equation}
Since $P_iI_i=I_i, \ \pi_iP_iI_i=\text{id}_{\C}$, we can rewrite (\ref{E:4})
as
$$
\eta_{h,i}\pi_ih_iP_iI_i=\lambda_i\eta_{h,i}\pi_iP_iI_i,\ \  {\text{or}}\ \
\eta_{h,i}\pi_ih_iP_i=\lambda_i\eta_{h,i}\pi_iP_i.
$$
But $h_iP_i=P_ih_i$, and so
\begin{equation}\label{E:5}
\eta_{h,i}\pi_iP_ih_i=\lambda_i\eta_{h,i}\pi_iP_i.
\end{equation}
The equation (\ref{E:5}), in turn, is equivalent to
\begin{equation}\label{E:6}
\eta_{h,i}\pi_i\xi_ip_ih=\lambda_i\eta_{h,i}\pi_i\xi_ip_i.
\end{equation}
We denote
\begin{equation}\label{E:7}
\theta_{h,i}=\eta_{h,i}\pi_i\xi_ip_i,
\end{equation}
a map from a neighborhood of $0\in\C^n$ into $\C$.
Then (\ref{E:6}) becomes $\theta_{h,i}h=\lambda_i\theta_{h,i}$. Now we define
$$
\theta_h=(\theta_{h,1},\dots, \theta_{h,n}),
$$
which is a germ of an analytic map at 0. This germ linearizes $h$:
$$
\theta_hh=(\theta_{h,1}h,\dots,\theta_{h,n}h)=(\lambda_1\theta_{h,1},\dots,
\lambda_n\theta_{h,n})=\Lambda \theta_h,
$$
where $\Lambda=\text{diag}(\lambda_1,\dots, \lambda_n)$ is an invertible
diagonal matrix, which has eigenvalues of $dh(0)$ on its diagonal.

The germ $\theta_h$ is biholomorphic. Indeed,
$$
\theta_{h,i}=\eta_{h,i}\pi_i\xi_ip_i=\eta_{h,i}\pi_iP_i\xi_i,\ \ i=1,\dots, n.
$$
Using the chain rule, we see that $d\theta_h(0)=A$, where $A$ is an
invertible diagonal matrix that diagonalizes $P_i$. We conclude that
$\theta_h$ is biholomorphic, and Lemma~\ref{L:Lin} is proved. $\Box$

\section{Simultaneous Linearization }\label{S:Sim}

Using Lemma~\ref{L:Lin}, we can linearize elements of $P_f(\Omega)$.
Namely, for every
$h\in P_f(\Omega)$ there exists $\theta_h$ (constructed in
Section~\ref{S:Sys}), such that $\theta_hh=\Lambda_h\theta_h$, where
$\Lambda_h$ is an invertible diagonal matrix. In particular, we can
linearize $f$:
$$
\theta_ff=\Lambda_f\theta_f,
$$
where the germ $\theta_f$ is biholomorphic at 0, and $\Lambda_f$ is an
invertible diagonal matrix.

\begin{lemma}\label{L:Slin}
For every $h\in P_f(\Omega)$ we have $\theta_h=\theta_f$.
\end{lemma}

\emph{Proof.}
Let us consider the germ
\begin{equation}\label{E:8}
\theta=\Lambda_f^{-1}\theta_hf,
\end{equation}
which is clearly biholomorphic. We have
$$
\theta h=\Lambda_f^{-1}\theta_hfh=\Lambda_f^{-1}\theta_hhf=\Lambda_f^{-1}
\Lambda_h\theta_hf=\Lambda_h\Lambda_f^{-1}\theta_hf=\Lambda_h\theta.
$$
Using (\ref{E:8}), we write the equation $\theta h=\Lambda_h\theta$ in the
coordinate form:
$$
(1/{\lambda_{f,i}})\theta_{h,i}fh=({\lambda_{h,i}}/{\lambda_{f,i}})
\theta_{h,i}f, \ \ i=1,\dots, n.
$$
By (\ref{E:7}) and the definition of $\xi_i$,
$$
(1/{\lambda_{f,i}})\eta_{h,i}\pi_iP_if_ih_i=({\lambda_{h,i}}/{\lambda_{f,i}})
\eta_{h,i}\pi_i P_if_i,\ \ i=1,\dots, n,
$$
where $f_i=\xi_if\xi_i^{-1}$.
Using the commutativity relations $f_iP_i=P_if_i,\ h_iP_i=P_ih_i$, which hold
since $\{p_i\}\subset S_f, \ h\in P_f(\Omega)$, we get
\begin{align}
&(1/{\lambda_{f,i}})\eta_{h,i}\pi_if_ih_iP_i=({\lambda_{h,i}}/{\lambda_{f,i}})
\eta_{h,i}\pi_if_iP_i,\ \ \text{or} \notag\\
&(1/{\lambda_{f,i}})\eta_{h,i}\pi_if_ih_iI_i=({\lambda_{h,i}}/{\lambda_{f,i}})
\eta_{h,i}\pi_if_iI_i, \ \ i=1,\dots, n.\notag
\end{align}
This is the same as
$$
((1/{\lambda_{f,i}})\eta_{h,i}\pi_if_iI_i)(\pi_ih_iI_i)=\lambda_{h,i}
((1/{\lambda_{f,i}})\eta_{h,i}\pi_if_iI_i), \ \ i=1,\dots, n,
$$
since $h_i$ locally preserves the $i$'th coordinate axis ($h_iP_i=P_ih_i$).
It is easily seen that
\begin{align}
&((1/{\lambda_{f,i}})\eta_{h,i}\pi_if_iI_i)(0)=0, \notag \\
&((1/{\lambda_{f,i}})\eta_{h,i}\pi_if_iI_i)'(0)=1.\notag
\end{align}
A normalized solution to a Schr\"oder equation is unique though. Thus we have
$$
\eta_{h,i}(\pi_if_iI_i)=\lambda_{f,i}\eta_{h,i}, \ \ \eta_{h,i}(0)=0,
\ \eta_{h,i}'(0)=1.
$$
Using the uniqueness argument again, we obtain $\eta_{h,i}=\eta_{f,i}$, and
hence $\theta_h=\theta_f$. The lemma is proved. $\Box$

According to Lemma~\ref{L:Slin}, the single biholomorphic germ $\theta_f$
conjugates the subsemigroup $P_f(\Omega)$ to some subsemigroup $D_f$ of invertible
diagonal matrices in $D_n$, the set of all $n\times n$ diagonal matrices
with entries in $\C$. We show that $D_f$ contains all invertible
diagonal matrices with sufficiently small entries. To do this, first we
extend $\theta_f$ to an analytic map on the whole domain $\Omega$ using
the formula
$$
\theta_f=\Lambda_f^{-l}\theta_ff^l,
$$
where $l$ is chosen so large that $\text{Cl}\{f^l(\Omega)\}$ is contained in
a neighborhood of 0 where $\theta_f$ is originally defined and biholomorphic;
the symbol Cl denotes closure. From the procedure of extending $\theta_f$ to
$\Omega$ we see that it is one-to-one and bounded in the whole domain.

Now, let $\Lambda=\text{diag}(\lambda_1,\dots,\lambda_n)$ be a matrix such
that $\text{Cl}\{\Lambda\theta_f(\Omega)\}\subset W$, where $W$ is a
neighborhood of $0\in\C^n$ for which $\text{Cl}\{\theta_f^{-1}W\}
\subset\Omega$. Such a matrix $\Lambda$ exists since $\theta_f$ is bounded
in $\Omega$. Consider $h=\theta_f^{-1}\Lambda\theta_f$, which belongs to
$G_0(\Omega)$. The map $h$ commutes with $f$ and all $p_i$'s. Indeed, using
the formula $\theta_ff\theta_f^{-1}=\Lambda_f$, we conclude that $hf=fh$ is
equivalent to $\Lambda\Lambda_f=\Lambda_f \Lambda$, which is a true relation
since both matrices $\Lambda$ and $\Lambda_f$ are diagonal. The relations
$hp_i=p_ih, \ i=1,\dots, n$, are verified similarly, using the formula
$\theta_fp_i\theta_f^{-1}=P_i$, which follows from the definition of
$\theta_f$.

\section{Solving a Matrix Equation}\label{S:Sol}

We proved that for an element $f\in VG_0(\Omega)$ there exists a
biholomorphic germ $\theta_f$ conjugating the semigroup $P_f(\Omega)$
to a subsemigroup $D_f\subset D_n$, which contains all invertible diagonal
matrices with sufficiently small entries.

\subsection{Conjugations $L$ and $R$}

Let $f\in VG_0(\Omega_1)$, and $g=\varphi f$. Then $g\in VG_0(\Omega_2)$,
and there is an isomorphism
$$
\Phi:\ \ S_f\rightarrow S_g.
$$
For the mappings $f$ and $g$ we have
$$
\theta_ff=\Lambda_f\theta_f,\ \ \theta_gg=M_g\theta_g,
$$
where $\Lambda_f, \ M_g$ are invertible diagonal matrices.

Let us consider the germ $L=\theta_g\psi \theta_f^{-1}$. This germ conjugates
the semigroups $D_f, \ D_g$:
\begin{align}
L\Lambda L^{-1}
&=\theta_g\psi\theta_f^{-1}\Lambda\theta_f\psi^{-1}
\theta_g^{-1} \notag\\
&=\theta_g\psi h\psi^{-1}\theta_g^{-1}=\theta_g j
\theta_g^{-1}=M,\notag
\end{align}
where $h\in P_f, \ \theta_fh=\Lambda\theta_f,\ j=\varphi h$, and $\theta_gj=M
\theta_g$.

Define $R(\Lambda)=L\Lambda L^{-1}$. Then $R:\ D_f\rightarrow D_g,$
$$
R(\Lambda_1\Lambda_2)=R(\Lambda_1)R(\Lambda_2), \ \ \Lambda_1, \ \Lambda_2
\in D_f.
$$

In what follows, we will identify $D_n$ with the multiplicative semigroup
$\C^n$ ($D_n\cong \C^n$) in the obvious way and consider a topology
on $D_n$ induced by the standard topology on $\C^n$.

\subsection{Extending $R$}

We are going to extend $R$ to an isomorphism of $D_n$. First,
we denote by $\overline{D}_f,\ \overline{D}_g$ the closures of $D_f,\ D_g$
in $D_n$, and for $\Lambda\in \overline{D}_f$ we set
$$
R(\Lambda)=\lim R(\Lambda_k),\ \ \Lambda_k\rightarrow\Lambda,\ \Lambda_k
\in D_f.
$$
This limit exists and does not depend on the sequence $\{\Lambda_k\}$, which
follows from the fact that $\psi^{\pm1}, \theta_f^{\pm1}$, and
$\theta_g^{\pm1}$
are continuous. The map $R$ is an isomorphism of topological semigroups
$\overline{D}_f$ and $\overline{D}_g$ (the inverse of $R$ has a similar
representation).

Next, we extend the map $R$ to $D_n$ as
$$
R(\Gamma)=R(\Gamma\Lambda)R(\Lambda)^{-1}, \ \ \Gamma\in D_n,
$$
where $\Lambda\in D_f$ is chosen so that $\Gamma\Lambda\in\overline{D}_f$.
This definition does not depend on the choice of $\Lambda$. Indeed, since
all matrices in question are diagonal (hence commute), the relation
$R(\Gamma\Lambda_1)R(\Lambda_1)^{-1}=R(\Gamma\Lambda_2)R(\Lambda_2)^{-1}$
is equivalent to $R(\Gamma\Lambda_1)R(\Lambda_2)=R(\Gamma\Lambda_2)
R(\Lambda_1)$, which holds.

The extended map $R$ is clearly an isomorphism of $D_n$ onto itself. Thus
we have
\begin{equation}\label{E:9}
R(\Lambda'\Lambda'')=R(\Lambda')R(\Lambda''), \ \ \Lambda', \ \Lambda''\in
D_n.
\end{equation}

Injectivity of $R$ and (\ref{E:9}) imply that $R(\Delta_i)=\Delta_j$ for
all $i$,
where $j=j(i)$ depends on $i$, and $j(i)$ is a permutation on $\{1,\dots, n\}$
(we recall that $\Delta_i=\text{diag}(0,\dots,1,\dots,0)$). This is because
$\{\Delta_i\}_{i=1}^n$ is the only system in $D_n$ with the following
relations: $\Delta_i\neq 0,\ \Delta_i^2=\Delta_i$, and $\Delta_i\Delta_j=0$
for $i\neq j$.

\subsection{A System of Scalar Equations}

Since all matrices $\Lambda$ and their images $R(\Lambda)$ are diagonal, we
can consider the matrix equation (\ref{E:9}) as $n$ scalar equations:
\begin{equation}\label{E:10}
r_j(\lambda_1'\lambda_1'',\dots, \lambda_n'\lambda_n'')=r_j(\lambda_1',
\dots, \lambda_n')r_j(\lambda_1'',\dots, \lambda_n''), \ \ j=1,\dots, n,
\end{equation}
where $r_j$ are components of $R$.
If we rewrite the equation $R(\Delta_i\Lambda)=\Delta_jR(\Lambda)$ in the
coordinate form, we see that $r_j(\lambda_1,\dots, \lambda_n)=r_j(0,\dots,
\lambda_i,\dots, 0)=q_j(\lambda_i)$; that is, each $r_j$ depends on only one
of the $\lambda_i$'s.
For each $j$ the corresponding equation in (\ref{E:10}) in terms of the
$q_j$'s
becomes
$$
q_j(\lambda_i'\lambda_i'')=q_j(\lambda_i')q_j(\lambda_i'').
$$
This equation has (\cite{aE93}, p. 130) either the constant solution
$q_j(\lambda_i)=1$, or
$$
q_j(\lambda_i)=\lambda_i^{\alpha_{ij}}\overline{\lambda}_i^{\beta_{ij}},
\ \ \alpha_{ij},\  \beta_{ij}\in \C, \ \ \alpha_{ij}-\beta_{ij}=\pm1.
$$

\subsection{Explicit Expression for $L$}

Going back to the function $L$, we have
\begin{equation}\label{E:11}
\begin{aligned}
L\text{diag}(\lambda_1,\dots, \lambda_n)=\text{diag}
(\lambda_{i(1)}^{\alpha_1}\overline{\lambda}_{i(1)}^{\beta_1},
\dots, \lambda_{i(n)}^{\alpha_n}\overline{\lambda}_{i(n)}^{\beta_n})L, \\
 \alpha_i-\beta_i=\pm1, \ \ i=1,\dots, n,\notag
\end{aligned} 
\end{equation}
where $i(j)$ is the inverse permutation to $j(i)$.

Let us choose and fix $(\mu_1,\dots, \mu_n)$ such that $(1/\mu_1,
\dots, 1/\mu_n)$ belongs to a neighborhood $W_0$ of $0\in\C^n$
where $L$ is defined, and let $W_1$ be a neighborhood of $0\in\C^n$
such that $(\mu_1z_1,\dots, \mu_nz_n)\in W_0$, whenever
$(z_1,\dots, z_n)\in W_1$.
Then from (\ref{E:11}) we have
\begin{align}
L&(z_1,\dots, z_n)
=L\text{diag}(\mu_1 z_1,\dots, \mu_n z_n)(1/\mu_1,\dots, 1/\mu_n) \notag\\
&=\text{diag}((\mu_{i(1)}z_{i(1)})^{\alpha_1}(\overline{\mu_{i(1)}
z_{i(1)}})^{\beta_1},\dots, (\mu_{i(n)}z_{i(n)})^{\alpha_n}
(\overline{\mu_{i(n)} z_{i(n)}})^{\beta_n}) \notag\\
&\times L (1/\mu_1,\dots, 1/\mu_n)=B (z_1^{\alpha_1}
\overline{z}_1^{\beta_1}, \dots, z_n^{\alpha_n}\overline{z}_n^{\beta_n}),
\notag
\end{align}
where $B$ is a constant matrix. The last formula is the explicit
expression for $L$.

\section{Proving that $\psi$ is (Anti-) Biholomorphic}\label{S:Pr}

To prove that $\psi$ is (anti-) biholomorphic is the same as to prove that
$L$ is (anti-) biholomorphic, because the relation $L=\theta_g\circ\psi\circ
\theta_f^{-1}$ holds. We showed that
\begin{equation}\label{E:12}
L(z_1,\dots, z_n)=B (z_1^{\alpha_1}\overline{z}_1^{\beta_1}, \dots,
z_n^{\alpha_n}\overline{z}_n^{\beta_n}), \ \  \alpha_i-\beta_i=\pm1,
\ \ i=1,\dots, n,
\end{equation}
in a neighborhood $W_1$ of $0$. From the representation (\ref{E:12}) we see
that $L$ is $\R$-differentiable and non-degenerate in $W_1\setminus
\cup_{k=1}^n
\{(z_1, \dots, z_n):\ z_k=0\}$. Since this is true for every point in the
domain $\Omega_1$, the map $\psi$ is $\R$-differentiable and
non-degenerate everywhere, with the possible exception of an analytic set.
Let us remove this set from $\Omega_1$, as well as its image under $\psi$
from $\Omega_2$. We call the domains obtained in this way $\Omega'$ and
$\Omega''$. Now the map $\psi:\ \Omega'\rightarrow\Omega''$ is
$\R$-differentiable and non-degenerate everywhere. It is clear that
if we prove that $\psi$ is (anti-) biholomorphic between $\Omega'$ and
$\Omega''$, then it is (anti-) biholomorphic between $\Omega_1$ and $\Omega_2$
due to a standard continuation argument \cite{sK82}. So we can think that
$\psi$ is $\R$-differentiable and non-degenerate in $\Omega_1$ itself.
The map $L$ thus has to be $\R$-differentiable and non-degenerate at 0.
However, this is the case if and only if $\alpha_i+\beta_i=1,
\ i=1,\dots, n$. Together with the equation $\alpha_i-\beta_i=\pm1$ it gives
us that either $\alpha_i=1,\ \beta_i=0$, or $\alpha_i=0, \ \beta_i=1$.

It remains to show that either $\alpha_i=1$ and $\beta_i=0$, or $\alpha_i=0$
and $\beta_i=1$, simultaneously for all $i$. Suppose, by way of contradiction,
that we have $L(z_1,\dots, z_n)=B (\dots, z_i,\dots, \overline{z}_j,\dots)$.
Then
$$
L^{-1}(w_1,\dots, w_n)=(\dots, l_i(w_1,\dots, w_n),\dots,
l_j(\overline{w}_1,\dots, \overline{w}_n),\dots),
$$
where $l_i,\ l_j$ are linear analytic functions. Let us look at an
endomorphism $f_0$ of $\Omega_1$ of the form
$$
f_0=\theta_f^{-1}\lambda(\dots, \theta_{f,i}\theta_{f,j},\dots,
\theta_{f,j},\dots)\theta_f,
$$
where $\theta_{f,i}\theta_{f,j}$ is in the $i$'th place, $\theta_{f,j}$
in the $j$'th, and $|\lambda|$ is sufficiently small.
Using (\ref{E:Con}) and the definition of $L$, we have
$$
\theta_g\varphi f_0\theta_g^{-1}=\theta_g\psi f_0\psi^{-1}\theta_g^{-1}=
L\theta_f f_0 \theta_f^{-1}L^{-1}.
$$
Thus,
\begin{align}
\theta_g\varphi f_0\theta_g^{-1}&(w_1,\dots, w_n) \notag\\
&= B' (\dots, l_i(w_1,\dots, w_n) l_j(\overline{w}_1,\dots, \overline{w}_n),
\dots, \overline{l}_j(w_1,\dots, w_n),\dots),\notag
\end{align}
for some constant matrix $B'$. This map, and hence $\varphi f_0$, is not
analytic though in a neighborhood of 0, which is a contradiction. Thus $L$,
and hence $\psi$, is either analytic or antianalytic in a neighborhood of 0.

Theorem~\ref{T:Mt} is proved completely. $\Box$

\section{Generalization}\label{S:Pr2}


Theorem~\ref{T:Mt} can be slightly generalized. Namely one may assume
that $\varphi$ is an epimorphism. We prove the
following theorem.

\begin{theorem}\label{T:Gt}
If $\varphi:\ E(\Omega_1)\rightarrow E(\Omega_2)$ is an epimorphism
between semigroups, where $\Omega_1,\ \Omega_2$ are bounded domains
in $\C^n,\ \C^m$ respectively, then $\varphi$ is an isomorphism.
\end{theorem}

{\emph{Proof}}
Since $\varphi$ is an epimorphism, it takes constant endomorphisms of
$\Omega_1$ to constant endomorphisms of $\Omega_2$, which follows
from~(\ref{E:Co}). Thus we can define a map $\psi:\ \Omega_1\to\Omega_2$
as in~(\ref{E:Dp}). Following the same steps as in verifying~(\ref{E:Con}),
we obtain
\begin{equation}\label{E:Sc}
\varphi f\circ\psi = \psi\circ f,\ \ {\text{for all }} f\in E(\Omega_1).
\end{equation}
We will show that~(\ref{E:Sc}) implies the bijectivity of $\psi$. The map
$\psi$ is surjective. Indeed, let $w\in\Omega_2$, and $c_w$ be the
corresponding constant endomorphism. Since $\varphi$ is an
epimorphism, there exists $f\in E(\Omega_1)$, such that
$\varphi f= c_w$. If we plug this $f$ into~(\ref{E:Sc}), we get
$$
\psi f(z) = w,
$$
for all $z\in \Omega_1$. Thus $\psi$ is surjective.

To prove that $\psi$ is injective, we show that for every
$w\in \Omega_2$, the full preimage of $w$ under $\psi$, $\psi^{-1}(w)$,
consists of one point.

Assume for contradiction that $S_w=\psi^{-1}(w)$ consists of more than
one point for some $w\in\Omega_2$. The set $S_w$ cannot be all of $\Omega_1$,
since $\psi$ is surjective. For $z_0\in\partial S_w\cap\Omega_1$, we
can find $z_1\in S_w$ and $\zeta\notin S_w$ which are arbitrarily close
to $z_0$. Let $z_2$ be a fixed point of $S_w$ different from $z_1$.
Consider a homothetic transformation $h$ such that $h(z_1)=z_1,\ h(z_2)=
\zeta$. Since the domain $\Omega_1$ is bounded, we can choose points
$z_1$ and $\zeta$ sufficiently close to each other so that $h$ belongs to
$E(\Omega_1)$. Applying~(\ref{E:Sc}) to $h$ we obtain
\begin{align}
&\varphi h(w)=\varphi h\circ \psi(z_1)=\psi\circ h(z_1)=\psi(z_1)=w, \notag\\
&\varphi h(w)=\varphi h\circ \psi(z_2)=\psi\circ h(z_2)=\psi(\zeta)\neq w.
\notag
\end{align}
The contradiction shows injectivity of $\psi$. Thus we have proved that
$\psi$ is bijective.

According to~(\ref{E:Sc}) we have
\begin{equation}\label{E:C1}
\varphi f=\psi\circ f\circ \psi^{-1}, \ \ \text{for all}\ f\in E(\Omega_1),
\notag
\end{equation}
which implies that $\varphi$ is an isomorphism.

Theorem~\ref{T:Gt} is proved. $\Box$




%
%
%

\bibliography{all}

\newpage
\vspace*{2cm}
\begin{center}
This page deliberately left blank
\end{center}

%
%
%

\appendix

\chapter{Speiser Graph is not Enough}
\label{S:Noten}

Here we give an example of a hyperbolic surface $(X, f)
\in F_3$, whose Speiser graph is parabolic. Thus it is essential
to consider the extended Speiser graph to determine the type of a
surface of the class $F_q$.

First, we consider a tree $D$, each vertex of which has degree 3.
Next, we fix a vertex $v\in VD$, and substitute each edge of $D$,
whose endpoints are at a distance $(n-1)$ and $n$ from $v$, by
$l_n$ edges in series, where $\{l_n\}$ is a sequence of odd
natural numbers. We complete  the graph obtained in this
way by edges, so that every vertex has degree 3, and there are no
algebraic elementary regions. Since all $l_n$
are odd, this is possible. The resulting graph is a Speiser graph
$\Gamma$.
We label the faces of $\Gamma$ by $0, 1, \infty$, and consider the
surface $(X, f)\in F(0, 1, \infty)$, corresponding to $\Gamma$
(the base curve is the extended real line).

If the sequence $\{l_n\}$ is increasing, then $(X, f)$
has a hyperbolic type, if and
only if (\cite{lT47}, \cite{lV50}, \cite{hW68})
\begin{equation}\label{E:Cons}
\sum_{n=1}^{\infty}\frac{\log{l_n}}{2^n}<\infty.
\end{equation}

On the other hand, using~(\ref{E:Exl}), we conclude that
$\lambda(T_v)\gtrsim\sum_{n=1}^{\infty}l_n/2^n$, where $T_v$ is
the family of paths in $\Gamma$ connecting $v$ to infinity.
Therefore, $\Gamma$ is parabolic if
\begin{equation}\label{E:Dse}
\sum_{n=1}^{\infty}\frac{l_n}{2^n}=\infty.
\end{equation}

We choose $l_n=2^n+1$. Combining (\ref{E:Cons}) and (\ref{E:Dse}),
we obtain a surface of a hyperbolic type, whose Speiser graph is
parabolic.

A straightforward computation shows that the mean excess of the
Speiser graph $\Gamma$ is 0.

\newpage
\vspace*{2cm}
\begin{center}
This page deliberately left blank
\end{center}

\chapter{Upper Mean Excess}\label{Upmeanex}

Here we show that the upper mean excess of every infinite Speiser
graph $\Gamma\in F_q$ is $\leq 0$. This proof is due to Byung-Geun Oh.

Let $\Gamma_{(i)}'$ be a double of $\Gamma_{(i)}$, i.e. a graph embedded
in a compact Riemann surface
obtained as follows. We take two copies of $\Gamma_{(i)}$,
one located above the other, and join by $s$ edges every pair of
boundary vertices of degree $q-s$ that are
located on the same vertical line. Let $n_f$ denotes the number of
edges on the boundary of a face $f\in F\Gamma_{(i)}'$.
A subgraph of $\Gamma_{(i)}'$, which is a copy of
$\Gamma_{(i)}$ (say a bottom copy) will again be denoted by $\Gamma_{(i)}$.
If a face
$f'\in F\Gamma_{(i)}'$ is induced by a face $f\in F\Gamma$,
then, clearly $n_{f'}/2\leq k_f$.
Therefore
\begin{equation}
\begin{aligned}
&2-\frac1{|V\Gamma_{(i)}|}\sum_{v\in V\Gamma_{(i)}} \sum_{\{f\in F\Gamma:
\ v\in
Vf\}}\bigg(1-\frac1{k_f}\bigg) \\
&\leq 2-\frac1{|V\Gamma_{(i)}|}\sum_{v\in
V\Gamma_{(i)}}\sum_{\{f\in F\Gamma_{(i)}':\ v\in Vf\}}
\bigg(1-\frac2{n_f}\bigg).
\end{aligned}
\end{equation}
If we assign the same value $\sum_{\{f\in F\Gamma_{(i)}':\ v\in Vf\}}
(1-2/{n_f})$ to every vertex $v'\in V\Gamma_{(i)}'$ that
lies above $v\in V\Gamma_{(i)}$, then (A.2)
is equal to
\begin{equation}\label{E:eq2}
\begin{aligned}
&2-\frac1{|V\Gamma_{(i)}'|}\sum_{v\in
V\Gamma_{(i)}'}\sum_{\{f\in F\Gamma_{(i)}':\ v\in Vf\}}
\bigg(1-\frac2{n_f}\bigg)\\
&=2-\frac1{|V\Gamma_{(i)}'|}\sum_{f\in F\Gamma_{(i)}'}\bigg(1-\frac2{n_f}
\bigg)n_f\\
&=2-\frac1{|V\Gamma_{(i)}'|}\sum_{f\in F\Gamma_{(i)}'}(n_f-2)
=2-\frac1{|V\Gamma_{(i)}'|}(2|E\Gamma_{(i)}'|-2|F\Gamma_{(i)}'|)\\
&\leq 2-\frac1{|V\Gamma_{(i)}'|}(2|V\Gamma_{(i)}'|-4)
=\frac4{|V\Gamma_{(i)}'|},
\end{aligned}
\end{equation}
where the inequality in (A.6) holds by the Euler
polyhedron formula ($|V| -|E| + |F|\leq 2$).
Since $|V\Gamma_{(i)}'|$ tends to infinity with $i$, the desired inequality
is established. $\Box$

\chapter{A Property of Extremal Length}\label{Proextlen}

Let $T,\ T_i,\ i\in I$, be families of paths in $G$, where $I$ is
at most countable. We assume that $ET_i\cap ET_j=\emptyset,\ i\neq j$.
Suppose that for every $t\in T$ and every $i\in I$, there exists
$t_i\in T_i$, which is a subpath of $t$. Then
\begin{equation}
\lambda(T)\geq\sum_{i\in I}\lambda(T_i).
\end{equation}

{\emph{Proof}} We can exclude from our consideration
the trivial cases when the sum on the right is zero, or when one
of the terms is infinite. For every $\epsilon>0$, and every $i\in
I$, we choose a density function $\mu_i$ on $ET_i$, such that for
every $t_i\in T_i$,
$$
\sum_{e\in Et_i}\mu_i(e)\geq 1,\ \ \sum_{e\in ET_i}\mu_i^2(e)
\leq\lambda(T_i)^{-1}+\epsilon.
$$
We choose a density function $\mu$ on $ET$, so that
$$
\mu(e)=\frac{\lambda(T_i)}{\sum_{j\in I}\lambda(T_j)}\mu_i(e),\ \ e\in ET_i,
$$
and 0 elsewhere.
Then for every $t\in T$,
$$
\sum_{e\in Et}\mu(e)=\sum_{i\in I}\sum_{e\in Et_i}
\frac{\lambda(T_i)}{\sum_{j\in I}\lambda(T_j)}\mu_i(e)
\geq \sum_{i\in I}\frac{\lambda(T_i)}{\sum_{j\in I}\lambda(T_j)}=1.
$$
Also,
$$
\begin{aligned}
&\sum_{e\in ET}\mu(e)^2=\sum_{i\in I}\sum_{e\in ET_i}\frac{\lambda(T_i)^2}
{(\sum_{j\in I}\lambda(T_j))^2}\mu_i(e)^2\\
&\leq\sum_{i\in I}
\frac{\lambda(T_i)^2}{(\sum_{j\in I}\lambda(T_j))^2}(\lambda(T_i)^{-1}+
\epsilon)
\leq \frac1{\sum_{i\in I}\lambda(T_i)} +\epsilon.
\end{aligned}
$$
Therefore,
$$
\lambda(T)^{-1}\leq\frac1{\sum_{i\in I}\lambda(T_i)} +\epsilon,
$$
and since $\epsilon$ is arbitrary, the desired inequality  is
established. $\Box$

\newpage
\vspace*{2cm}
\begin{center}
This page deliberately left blank
\end{center}


%
%
%

\begin{vita}
  Sergiy Merenkov was born on February 5, 1974, in
Cherepovets, Vologda region, USSR (present Russia). He 
obtained his M.S. degree in 1996, from 
Kharkov State University, Ukraine. In 1999, he was 
accepted to Purdue University.   
\end{vita}

\end{document}
